
%
%

 \documentstyle{amsppt}
\pagewidth{6.4in}\vsize8.5in\parindent=6mm
\parskip=3pt\baselineskip=14pt\tolerance=10000\hbadness=500
\NoRunningHeads
\loadbold
\topmatter
\title
Singular maximal functions and Radon transforms 
near
 $\text{L}^{\text 1}$
\endtitle
\author
Andreas Seeger \ \ \  Terence Tao \ \ \ James Wright
\endauthor
\date 
May 12, 2002
\enddate
\abstract
We show that some singular maximal functions and singular Radon transforms
 satisfy a weak type $L\log\log L$ inequality.
Examples include the maximal function and Hilbert 
transform associated to averages
along a parabola. 
The weak type inequality  yields pointwise convergence results
for functions which are locally in
$L\log\log L$. 
\endabstract
\thanks
The first author is supported in part by a grant from the National Science
Foundation.  The second author is a Clay Prize fellow
and is supported by the Sloan and Packard foundations.
  \endthanks
\address
Department of Mathematics,
University of Wisconsin, Madison, WI 53706-1388, USA
\endaddress
\email seeger\@math.wisc.edu\endemail
\address
Department of Mathematics,
University of California, Los Angeles, CA  90095-1555, USA
\endaddress
\email tao\@math.ucla.edu\endemail
\address
Department of Mathematics and Statistics,
University of Edinburgh,
King's Buildings,
Mayfield Rd.,
Edinburgh EH3 9JZ, U.K.
\endaddress
\email wright\@maths.ed.ac.uk\endemail
\endtopmatter
\document

\def\mx{{\max}}

\def\vth{\vartheta}

\def\fnkw{{f^{n,\ka}_w}}

\def\R{\Bbb R}
\def\Q{\Cal Q}

\def\emph#1{{\it #1 }}

\define\ga{\gamma}
\define\Th{\Theta}

\define\sgn{{\text{\rm sign }}}

\define\loc{{\text{\rm loc}}}


\define\dist{{\text{\rm dist}}}

\define\supp{{\text{\rm supp }}}

\define\inn#1#2{\langle#1,#2\rangle}
\define\biginn#1#2{\big\langle#1,#2\big\rangle}

\define\lcontr{\rfloor}
\define\lco#1#2{{#1}\lcontr{#2}}
\define\lcoi#1#2{\imath({#1}){#2}}
\define\rco#1#2{{#1}\rcontr{#2}}
\redefine\exp{{\text{\rm exp}}}
\define\bin#1#2{{\pmatrix {#1}\\{#2}\endpmatrix}}
\define\meas{{\text{\rm meas}}}

\define\card{\text{\rm card}}
\define\lc{\lesssim}


\define\del{\delta}             
\define\eps{\varepsilon}

\define\ka{\kappa}
            
\define\la{\lambda}             \define\La{\Lambda}

\define\fQ{{\frak Q}}
\define\fR{{\frak R}}
\define\fS{{\frak S}}

\define\fW{{\frak W}}


\define\fv{{\frak v}}
\define\fw{{\frak w}}

\define\bbR{{\Bbb R}}

\define\bbZ{{\Bbb Z}}

\define\cE{{\Cal E}}

\define\cG{{\Cal G}}
\define\cH{{\Cal H}}

\define\cJ{{\Cal J}}

\define\cM{{\Cal M}}

\define\cO{{\Cal O}}
\define\cP{{\Cal P}}
\define\cQ{{\Cal Q}}




\head{\bf 1.Introduction}
\endhead

Let $\Sigma$ be a compact  smooth   hypersurface of $\Bbb R^d$, 
and let $\mu$ be a compactly supported smooth
density on $\Sigma$, {\it i.e.} $$\mu=\chi d\sigma$$
where $\chi\in C^\infty_0(\Bbb R^d)$ and $d\sigma$ is the surface carried measure on $\Sigma$.

Unless stated otherwise we shall always make the following 
\proclaim{Curvature Assumption}
The  Gaussian curvature 
 does not vanish to infinite order on $\Sigma$.
\endproclaim

We consider a group
of dilations on $\Bbb R^d$, given by $t^P= \exp (P\log t)$, $t>0$,   and
 we assume that $P$ is a $d\times d$
matrix whose eigenvalues have positive real part. For $k\in \Bbb Z$
we set
$\del_k =2^{kP}$ and define the measure $\mu_k$ by
$$\inn{\mu_k}{f} =\inn{\mu}{f(\del_{k}\cdot)}.
\tag 1.1
$$
We shall consider the convolutions $\mu_k*f$ and  study the behavior of the maximal function
$$\cM f(x)=\sup_{k\in \Bbb Z} | \mu_k * f(x)|
\tag 1.2$$
and some related singular integrals. By a rescaling we may assume that the measure $\mu$ is supported in the unit ball $\{x:|x|\le 1\}$.

The first complete $L^p$ bounds ($1<p<\infty$)  for a class of such operators (Hilbert transforms on curves) seems to be due to Nagel, Rivi\`ere
and Wainger \cite{9}. A classical reference is the article by 
Stein and Wainger \cite{17} containing many related results; see also the paper by Duoandikoetxea and Rubio de Francia \cite{6} which contains 
general results for maximal functions and singular integrals 
generated by singular measures, with decay assumptions on the 
Fourier transform. Concerning the behavior on $L^1$ 
it is presently not known  even for the 
special classes considered 
here 
whether the maximal operator $\cM$ is of weak type $(1,1)$, {\it i.e.} 
whether it maps $L^1$ to the Lorentz space $L^{1,\infty}$. 
This question had been raised in \cite{17}.
For some 'flat' cases counterexamples are in \cite{3}, but these do not seem to apply 
in the case of our  curvature assumption.

We shall examine the behavior of the maximal function on 
spaces ``near'' $L^1$.
Two results  in this direction are known:
Christ and Stein \cite{4} showed by an extrapolation argument
that if $f$ is supported in a cube $Q$ and $f\in L\log L(Q)$ then the maximal function
$\cM f$ belongs to  $L^{1,\infty}$ (again under substantially weaker 
finite type assumptions).
Moreover Christ \cite{2} showed that the lacunary 
spherical maximal function
maps the standard  Hardy space $H^1(\bbR^d)$ to $L^{1,\infty}$, and 
that  maximal functions and Hilbert transforms 
associated to a parabola in $\Bbb R^2$
 map the appropriate Hardy space with respect to nonisotropic dilations to $L^{1,\infty}$.
Weak $L^1$ (see also Grafakos \cite{8} and our recent paper \cite{12} for related results). 
For the two operators associated to the  parabola  $(t,t^2)$ it is also known (\cite{11}) 
that they map the smaller product-type Hardy space
$H^1_{\text{prd}}(\bbR\times\bbR)$ to the smaller  Lorentz space $L^{1,2}$.

We recall that for $f$ to belong to a Hardy space $H^1$
a  rather substantial cancellation condition has to be satisfied.
If locally the cancellation is missing
one has a  restriction on the size of $f$; more  precisely
 if a function $f\in H^1$ is single signed in an open  ball
then $f$ belongs to $L\log L(K) $ for all  compact subsets $K$ of this ball. This can be deduced
 from the maximal function characterization of $H^1$ and  the
fact that $f_0\in L\log L(q_0)$ if $f_0$ is supported on the cube 
$q_0$ and
the appropriate variant of the Hardy-Littlewood maximal function of $f_0$
belongs to $L^1(q_0)$, see
\cite{15, \S I.5.2 (c)}.
Here we are  interested in the 
behavior in Orlicz spaces near $L^1$ without assuming
additional cancellation conditions.

Our  main result is that the maximal operator 
acts well on $L\log \log L$ and the global version 
satisfies weak type $L\log\log L$ inequalities. We first give a 

\definition{Definition} 
{\it
Let  $\Phi:\Bbb R^+\to \Bbb R^+$ be a  convex function and let
 $T$ be an operator mapping simple functions on $\Bbb R^d$  to measurable functions.
$T$ {\it is of weak type $\Phi(L)$} if  there is a constant $C$ so that the inequality
$$
\big|\{x\in \Bbb R^d: |T f(x)|>\alpha\}\big|
\le
\int \Phi\big(\frac{C|f(x)|}{\alpha}\big) dx
\tag 1.3
$$
holds for all $\alpha>0$.
\enddefinition

Abusing the notation slightly we shall say that $T$ 
 is {\it of weak type $L\log\log L$} if  there is a 
constant $C$ so that the inequality (1.3) holds with
$\Phi(t)= t\log\log(e^2+t).$  


\proclaim{Theorem 1.1}
The maximal operator $\cM$ is of weak type $L\log\log L$.
\endproclaim

We  also prove a related theorem on singular convolution operators 
with kernels supported on hypersurfaces (assuming our finite type curvature assumption).

Let $\mu_k$ be as in (1.1)  and   assume that in addition
$$\int d\mu=0.\tag 1.4 $$
For Schwartz functions $f$ define
 the singular integral operator (or singular Radon transform) $T$ by
$$T f(x)= \sum_{k\in \Bbb Z} \mu_k * f.
\tag 1.5$$ 
\proclaim{Theorem 1.2}
  $T$ extends to an 
operator which is of weak type
$L\log\log L$.
\endproclaim

\bigskip 

\subheading{\bf 1.3 Remarks and examples}

{\bf 1.3.1.} Theorem 1.1 implies an estimate on the Orlicz space
$\Phi(L)(Q_0)$ where $Q_0$ is a unit cube and
the norm on $\Phi(L)$ is given by
$\|f\|_{\Phi(L)}=\inf\{\alpha> 0: \int_{Q_0} \Phi(|f(x)|/\alpha) dx\le 1\}$.
Consider the local maximal operator
$${\cM}_\loc f(x)= \sup_{k<C}|\mu_k* [f\chi_{Q_0}](x)|;$$
then $\cM_\loc$ maps $L\log\log L(Q_0)$ to $L^{1,\infty}$.
To see this we may assume that
$\|f\|_{L\log\log L(Q_0)}=1$. Then the estimate
$$|\{x\in Q_0: {\cM_\loc}f>\alpha\}|\lc\alpha^{-1}$$
is trivial for $\alpha<1$ while for $\alpha>1$ it follows from the better estimate (1.3).

We note  that conversely the better estimate
$|\{x\in \bbR^n:\cM_\loc f>\alpha\}|\lc\int\Phi(C|f(x)|/\alpha)$
can be deduced from 
the $L\log\log L(Q_0)\to L^{1,\infty}$ boundedness
by the Orlicz space variant of Stein's theorem \cite{14}. 
Then the global variant of Theorem  1.1 follows by 
scaling and limiting arguments.

{\bf 1.3.2.} Similarly if we assume the cancellation condition (1.4) then 
the local singular Radon transform
$\sum_{k<C}\mu_k* [f\chi_{Q_0}](x)$ 
maps 
$L\log\log L(Q_0)$ to $L^{1,\infty}$.

{\bf 1.3.3.}  Suppose that $\int d\mu=1$ and suppose that the 
 measurable function $f$ belongs locally to
$L\log \log L$; i.e.
$\int_K |f(x)|\log\log (e^2+|f(x)|) dx<\infty$ for every compact set $K$. Then 
$\lim_{k\to-\infty}\mu_k*f(x)=f(x)$ almost everywhere.

This follows by a standard argument.
Observe that  we have
 $\int\alpha^{-1} |f(x)|\log\log(e^2+\alpha^{-1} |f(x)|) dx<\infty$, for every $\alpha>0$.
Fix  $\alpha>0$ and let
$$\Omega_{\alpha}(f)=\big\{x:\limsup_{k\to-\infty} \mu_k*f(x)-
\liminf_{k\to -\infty} \mu_k* f(x)>\alpha\big\}.$$
Given $\eps>0$ we  show that $|\Omega_{\alpha}(f)|<\eps$. One can find a bounded
 function $h$ with compact support
so that $\int\Phi(2C|f-h|/\alpha) dx \le \eps$ and  since 
$\mu_k* h\to h$ almost everywhere
we see that
$\Omega_{\alpha/2}(h)$ has measure zero.
Moreover $|\Omega_\alpha(f)|\le |\Omega_{\alpha/2}(f-h)|+
|\Omega_{\alpha/2}(h)|$ and by Theorem 1.1 
 we see that $\Omega_{\alpha/2}(f-h)$ and thus $\Omega_\alpha(f)$
  has measure $<2\eps$. Since $\eps$ was arbitrary we see that $\Omega_\alpha(f)$
has measure zero; thus  $\cup_{m}\Omega_{2^{-m}}(f)$ has measure zero
and the result on pointwise convergence follows.

{\bf 1.3.4.} Examples of Theorem 1.1 include the lacunary 
spherical maximal operator  where $\mu_k* f$ is 
the  average of $f$ 
over the sphere of radius $2^k$ centered at $x$
(for the early $L^p$ results see \cite{1}, \cite{5}). 
The sphere may
be replaced by any smooth compact hypersurface for which the  curvature 
vanishes of finite order only, and the isotropic dilations 
may be replaced by nonisotropic ones.
We remark that the proof of Theorem 1.1 for isotropic dilations 
is much less technical, see the expository note \cite{13}.

{\bf 1.3.5.} Other examples of Theorem 1.1 concern the 
averages along a parabola 
$$\cP_r f(x) =\frac 1 r\int_0^r f(x_1-t, x_2-t^b)dt$$
or higher dimensional versions for paraboloids $(t',|t'|^b)$, $b\neq 1$.
Again if  $f$ belongs locally to
$L\log \log L$ then $\lim_{r\to 0}\cP_r f(x)=f(x)$ almost everywhere.

{\bf 1.3.6.} Similarly Theorem 1.2 can be used to deduce the weak type 
$L\log\log L$ inequality for the Hilbert transform
$$\cH f(x) = p.v. \int_{-\infty}^\infty f(x_1-t, x_2-t^b)\frac{dt}t.$$


We give a brief outline of the paper.
The main novelty in this paper is a stopping time argument based on
the quantities of {\it thickness} $\Theta_n$ and {\it length} $\La_n$
 associated to a density $v(x) dx$ (depending on an additional parameter $n$). 
  Basically, the point is that the length $\La_n[v]$ is used to control the size of an
exceptional set
while the thickness $\Theta_n[v]$ is used to  control the $L^2$ norm of 
an essential part
of the maximal function outside of the exceptional set, for suitable choices of $v$.
The quantities of length and thickness are complementary in some sense; 
this and other  basic properties are   discussed in \S2. In
\S3 we include  preliminary and standard  arguments from
Calder\'on Zygmund theory. These arguments can be skipped by the experts; they may be used to 
 reprove the known   $L\log L$ estimates. 
In \S4 we describe the stopping time argument based on length and thickness.
The proof of the
weak-type  $L\log\log L$ inequality  for the maximal operator is given in \S5.
The bounds for the singular Radon transforms are discussed in \S6.

\head{\bf 2. Length and thickness}\endhead

In this section  let $v$ be  an integrable nonnegative function  which 
vanishes in the complement of a dyadic cube $q$.
 Dyadic cubes are supposed to be `half-open', i.e. of the form
$\prod_{i=1}^d [ n_i2^{m}, (n_i+1)2^m)$ where $n_i,m \in \bbZ$.

We
define a dyadic version of a 
{\it one-dimensional Hausdorff content}
or simply  {\it length } $\la(E)$
  to be
$$ \la(E) := \inf_{\cQ} \sum_{Q \in \cQ} l(Q)\tag 2.1$$
where $\cQ$ ranges over all finite collections $\cQ$ of dyadic cubes  
with
$E\subset\bigcup_{Q\in\cQ}Q$,
 and $l(Q)$ denotes the sidelength of $Q$.
Note that this definition differs from the usual definition of a one-dimensional
Hausdorff measure as $\la(E)\le l(Q)$ if $E$ is contained 
in the dyadic cube $Q$.

Given $n\in \Bbb Z$  we denote by $E_n[v]$ the conditional 
expectation of 
$v$,
for the $\sigma$-algebra generated by dyadic cubes of sidelength $2^{-n}$; thus
$$E_n[v](x)=\sum_Q \chi_Q(x) |Q|^{-1}\int_Q v(y) dy $$ where of course 
the sum runs over all  dyadic cubes of sidelength $2^{-n}.$  We  also define
$$\fS_n(v)=\{x: E_n[v](x)\neq 0\}.\tag 2.2$$
Notice that $v(x)=0$ for almost every $x\in \bbR^d\setminus \fS_n[v]$ 
since $v$ is nonnegative. Now define
$$\Lambda_n[v]=\lambda(\fS_n(v)).\tag 2.3$$
Note that $\fS_n(v)$ is a union of dyadic cubes of length $2^{-n}$ and 
therefore the infimum in the definition of $\lambda$ becomes a minimum; 
{\it i.e.} there is a collection $\cQ$ of dyadic cubes covering the set 
$\fS_n(v)$ so that
$ \La_n[v]= \sum_{Q \in \cQ} l(Q)$.
Here the cubes in $\cQ$ have to be chosen to be of sidelength at least $2^{-n}$.

Next 
we define the {\it thickness} of $v$
 to be the quantity
$$ \Th_n[v] := \sup_{Q} \frac{1}{l(Q)}\int_Q v(x) dx \tag 2.4$$
where $Q$ ranges over all dyadic cubes of sidelength $l(Q)\ge 2^{-n}$.
Clearly, if $v$ vanishes off  a dyadic cube $q$  it is sufficient to
only consider  dyadic subcubes of $q$ in (2.4).

We note that the restriction to dyadic cubes in the definition of length and thickness
is convenient but not essential. Since every cube of sidelength $2^L$ ($L\in \bbZ$) 
is contained in a union of $2^d$ dyadic cubes of sidelength $2^L$ we observe that
$$\aligned
&\Th_n[v(\cdot+a)]\le 2^d\Th_n[v]
\\
&\La_n[v(\cdot+a)]\le 2^d\La_n[v].
\endaligned
\tag 2.5
$$

The quantities of length and thickness are complementary.  
Namely, it is immediate from the definitions of $\La_n$ and $\Theta_n$
that 
$$ \int v(x) dx  \leq \La_n[v]\Th_n[v].
\tag 2.6$$
The bound (2.6) can be attained, for instance if
 $v$ is the characteristic function of  a dyadic box. It would be
 desirable to have a converse  to (2.6), with bounded constants,
but this generally does not hold as the following example shows.
Let $E_n$ be the union of $n+1$ rectangles $R_\nu$, parallel to the coordinate axes,
 with dimensions 
$(2^{-\nu},1)$ so 
that the left lower endpoint of $R_\nu$ has coordinates $(\nu,0)$, $\nu=0,\dots, n$. 
Let $v_n=\chi_{E_n}$.
Then $\La_n[v_n]=n+1$, $\int v_n(x) dx< 2$ and  $\Th_n[v_n]=1$; thus the converse of
(2.3) fails with a uniform constant. 

However we shall show that $v$ can be efficiently decomposed into 
a sum of functions for which a converse of (2.6) does  hold. 
The main result needed to achieve this is

\proclaim{Proposition 2.1} Let $q$ be a dyadic cube with $l(q)\ge 2^{-n}$.
Suppose that $v$ is a bounded nonnegative 
measurable function supported in  $q$.
Then there exists a decomposition 
$$v=g+h$$ with 
nonnegative functions  $g$ and $h$ and $g$, $h$  vanish in the complement of the set 
$\fS_n(v)\subset q$; moreover the inequalities 
$$\Lambda_n[h]\le \frac 12 \Lambda_n[v]
\tag 2.7$$
and
$$\Lambda_n[v] \Th_n[g] \leq 8\int g(x) dx \tag 2.8$$
hold.
\endproclaim

In particular we see from (2.7/8) that  the function $g$ satisfies 
$$\Lambda_n[g] \Th_n[g] \leq 8\int g(x) dx,$$ thus a converse to (2.6).

We shall first prove a technical result which states that for each dyadic cube one may construct
a function $v_I$ from $v$  so that $v_I$ has `controlled' thickness and `large' integral.

\proclaim{Lemma 2.2} Let $\ga>0$. 
For any dyadic cube $I$ of sidelength $\ge 2^{-n}$, 
there exists a (possibly empty) collection
$\Q[I]$ of disjoint dyadic cubes
of sidelength 
$\ge 2^{-n}$  contained in $I$,  and
a measurable function $v_I$ 
 such that
$$
0\le v_I(x)\le v\chi_I(x)\tag 2.9
$$
for all $x\in \Bbb R^d$, 
$$\Th_n[v_I] \leq 2\ga \tag 2.10$$
and
$$ 2\int v_I(x) dx \ge 2\ga \sum_{Q \in \Q[I]} l(Q) + 
\int_{I \backslash
\bigcup_{Q \in \Q[I]} Q} v(x) dx. \tag 2.11$$
\endproclaim

\demo{\bf Proof }
We prove this by induction on the sidelength of $I$.  
We first assume that $l(I)= 2^{-n}$. 
Notice that in this case  we have 
$$\Theta_n [v\chi_I]= 
\frac{1}{l(I)}\int_I v(x) dx.$$

We distinguish two cases. First
 if $\Theta_n[v\chi_I]\le 2\ga$ we choose $v_I=v\chi_I$ and take 
for $\Q[I]$ the empty collection. Clearly (2.9), (2.10), (2.11) are satisfied.

Next  if $\Theta_n[v\chi_I]> 2\ga$
 we may choose a measurable function $v_I$ which vanishes outside $I$ 
such that $0\le v_I(x)\le  v\chi_I(x)$ for all $x\in \bbR^d$ and 
$$\ga\le \frac{1}{l(I)}\int_I v_I(x) dx \le 2\ga.
\tag 2.12
$$ Clearly
$\Theta_n[v_I]\le 2\ga$. For $\Q[I]$ we take the singleton collection $\{I\}$ and (2.11)
 is satisfied 
because of the first inequality in (2.12).

Now 
fix a dyadic cube $I$ with $l(I) > 2^{-n}$ and suppose  that the lemma has been proven for all
proper dyadic subcubes
$I'$ of sidelength at least $2^{-n}$.  Partition $I$ into $2^d$ 
subcubes $I_1, \ldots, I_{2^d}$ of
sidelength $\frac{1}{2} l(I)$.  By the induction hypo\-the\-sis,
we may construct collections $\Q[I_j]$ and measurable functions  $v_{I_j}$ for $j=1, \ldots, 2^d$
satisfying the properties of the lemma relative to $I_j$.

To prove the assertion for $I$ we again distinguish two cases.
First suppose that
$$\sum_{j=1}^{2^d}\int v_{I_j}(x)dx \leq 2\ga l(I).
\tag 2.13$$
In this case we simply define $v_I(x):=\sum_{j=1}^{2^d} v_{I_j}(x)$ and
$\Q[I] := \bigcup_{j=1}^{2^d} \Q[I_j]$. Then  by the induction hypothesis
$$
2\int v_I(x) dx=\sum_{j=1}^{2^d} 2\int v_{I_j}(x) dx\ge
\sum_{j=1}^{2^d} \Big[2\ga \sum_{Q\in\cQ[I_j]}l(Q)+
\int_{ I_j \setminus \cup_{Q\in \cQ[I_j]}Q} v(x) dx\Big]
$$ which is equal to the right hand side of (2.11). From  (2.13) it follows that
$$\frac{1}{l(I)}\int v_I(x) dx \le 2\ga$$ 
and if $Q$ is a proper dyadic subcube of $I$ then $Q\subset I_j$ for
some $j$ and
$$\frac{1}{l(Q)}\int_Q v_I(x) dx=
\frac{1}{l(Q)}\int_Q v_{I_j}(x) dx \le 2\ga$$
 by the induction hypothesis. Altogether (2.10)
follows in case (2.13).

Now suppose that
$$\sum_{j=1}^{2^d}\int v_{I_j}(x)dx > 2\ga l(I).\tag 2.14$$
In this case we can find a function $v_I$ so that 
$v_I(x)\le \sum_{j=1}^{2^d} v_{I_j}(x)$ and
$$\ga l(I)\le \int v_{I}dx \le 2\ga l(I).\tag 2.15$$
We then take for $\cQ[I]$ the singleton set $\{I\}$. Then (2.11) is immediate by (2.15). 
Clearly also by (2.15) 
$\frac{1}{l(I)}\int v_I(x) dx \le 2\ga$. As above we can use the 
induction hypothesis to see 
that if $Q$ is a proper dyadic subcube, 
thus  contained in an $I_j$, we have
$\frac{1}{l(Q)}\int_Q v_I(x) dx\le 
\frac{1}{l(Q)}\int_Q v_{I_j}(x) dx \le 2\ga$, thus altogether (2.10) 
also holds  in this case.\qed
\enddemo

\demo{\bf Proof of Proposition 2.1}
We define the {\it critical thickness} $\vth_n(v)$
to be the
largest non-negative number $\ga$ such that the inequality
$$  \ga\Lambda_n[v] \leq 2\ga\sum_{Q \in \cQ} l(Q) +
\int_{q \backslash \bigcup_{Q \in \cQ} Q}v(x) dx \tag 2.16$$
holds for all finite collections $\cQ$ of dyadic cubes of sidelength $2^{-n}$ (here the empty collection is admitted).
Equivalently, one can define $\vth_n(v) $ by
$$ \vth_n(v)  := \inf_{\cQ} \frac{ \int_{q \backslash 
\bigcup_{Q \in \cQ} Q }v(x)dx}
{(\Lambda_n[v] - 2\sum_{Q \in \cQ} l(Q))_+}.\tag 2.17$$
Observe that since $v$ vanishes in the complement of 
$q$
and since all cubes have sidelength at least $2^{-n}$ 
we are in effect taking the infimum over a finite set of collections, each consisting of a finite number of cubes,
 so that this infimum becomes a minimum, and (2.16) holds with $\gamma=\vth_n(v)$.

Clearly $\vth_n(v)\le \Lambda_n[v]^{-1}\int v(x) dx$. Observe also  that  $\vth_n(v)  > 0$
since 
$\int_{q \backslash \bigcup_{Q \in \cQ} Q}v(x) dx $ is positive
whenever $\sum_{Q \in \cQ} l(Q) \leq \Lambda_n[v]/2$.

We can now find a  finite collection $\cQ_1$ of dyadic cubes in $q$, of sidelength at least $2^{-n}$, 
so that

$$  \vth_n(v) \Lambda_n[v] = 2\vth_n(v)  \sum_{Q \in \cQ_1} l(Q) + \int_{E_*} v(x) dx \tag 2.18$$
where 
$$ E_* := q \backslash \bigcup_{Q \in \Q_1} Q.
\tag 2.19
$$
We claim that
$$ \Th_n[v\chi_{E_*}] \leq 2\vth_n(v) .\tag 2.20$$
Indeed, suppose for contradiction that there existed a dyadic cube $Q'$ such that
$$ \int_{E_*\cap Q'}v(x) dx  >2 \vth_n(v)  l(Q').\tag 2.21$$
By (2.21) and $\vth_n(v)>0$ we have 
$|E_*\cap Q'|>0$ which implies that $Q'\notin\cQ_1$.
If we apply (2.16)
to the collection $\Q_1
\cup \{Q'\}$ we obtain
$$  \vth_n(v) \Lambda_n[v] \leq 2\vth_n(v) 
\Big(l(Q') + \sum_{Q \in \Q_1} l(Q)\Big) + \int_{E_*\setminus Q'} v(x) dx,
$$
but by (2.18) this implies 
$$\int_{E_*} v(x) dx\le 2\vth_n(v) l(Q')+\int_{E_*\setminus Q'} v(x) dx$$
contradicting  (2.21).  This proves (2.20).

We shall now invoke Lemma 2.2  with $\gamma=\vth_n(v)$ and $I=q$, thus 
finding a function  $v_q$ and a collection 
$\Q[q]$ obeying the properties in the lemma.  We define
$$
g(x) = v(x)\chi_{E_*}(x)+ v_q(x) \chi_{q\setminus E_*}(x)
$$
and
$$
h(x)=\big(v(x)-v_q(x)\big) \chi_{q\setminus E_*}(x).
$$
Observe that $g$ and $h$ are nonnegative functions.
To show (2.7) we use that
$\La_n[h]\le \la (q\setminus E_*)$ since the latter set is a union of 
dyadic cubes of sidelength $2^{-n}$. Thus we observe
$$ \La_n[h]\le \sum_{Q\in \cQ_1}l(Q)\le\frac 12 \Lambda_n[v],$$
by (2.18). This gives (2.7).

To show (2.8) we use that $v_q\le v$ and observe that by (2.11)
$$\int g(x) dx\ge \int v_q(x) dx\ge
 \frac{1}{2} \Big(2\vth_n(v) \sum_{Q \in \Q[q]} l(Q) + 
\int_{q\setminus \cup_{Q\in \cQ[q]} Q}  v(x) dx\Big ),$$
 since now $\ga=\vth_n(v)$.  By (2.16) we thus see that
$$\int g(x) dx \geq \frac 12 \La_n[v]\vth_n(v).$$
By (2.20) and (2.10) 
$$ \Th_n[g] \leq \Th_n[v\chi_{E_*}] + \Th_n[v_q] \leq 2\vth_n(v)  + 2\vth_n(v)  = 4\vth_n(v) ,$$
we see that $\Th_n[g]\le 8\La_n[v]^{-1}\int g(x) dx$ which is (2.8).\qed
\enddemo

\remark{\bf Remark} There are analogues of Proposition 2.1 where 
for $0<\beta<d$ the
length $\la(E)$ is replaced by the $\beta$-dimensional 
Hausdorff content $$\la_\beta(E)
= \inf_{\cQ} \sum_{Q \in \cQ} l(Q)^\beta$$
where again  $\cQ$ ranges over all finite collections $\cQ$ of 
dyadic cubes  
with
$E\subset\cup_{Q\in\cQ}Q$. Then if we define $\La_{\beta,n}(v)=\la_\beta(\fS_n(v))$ and the $\beta$-thickness by
$$ \Th_{\beta,n}[v] := \sup_{Q} \frac{1}{l(Q)^\beta}\int_Q v(x) dx 
$$ then an assertion analogous to Proposition 2.1 holds true. The proof requires only notational changes.

\endremark

In what follows it will be convenient to extend the definition of 
length and thickness to not necessarily nonnegative functions, and we simply put
$$\Lambda_n[f]:=\Lambda_n[|f|], \quad
\Th_n[f]:=\Th[|f|].
$$

Proposition 2.1 can be applied iteratively. This leads to
\proclaim{Proposition 2.3} Suppose that $f$ is  integrable and vanishes in the 
complement of dyadic cube of length $1$.
Set $h_0(x)=f(x)$.  For $m\ge 1$  we may decompose 
$$f=h_m+\sum_{\nu=1}^m g_\nu$$
almost everywhere, so that the following properties hold.

(i) $h_m(x)$ and the $g_\nu(x)$ are nonnegative if and only 
if $f$ is nonnegative, and 
$h_m(x)$ and the $g_\nu(x)$ are nonpositive if and only 
if $f$ is nonpositive.


(ii) 
$\Theta_n[g_\nu]\La_n[h_{\nu-1}]\le 8\int |g_\nu(x)|dx$.

(iii) $\Lambda_n[h_m]\le 2^{-m}\Lambda_n[f]$.

(iv) If  $m\ge n$ then $g_{m+1}=h_m$, $h_{m+1}=0$.

\endproclaim

\demo{\bf Proof} 
%
We first extend the statement of Proposition 2.1 to not necessarily 
nonnegative functions, in the obvious way. We simply decompose 
$|f|=\tilde g+\tilde h$ as in Proposition 2.1, and then define
$g(x)=\tilde g(x)\sgn\!(f(x))$, and
$h(x)=\tilde h(x)\sgn\!(f(x))$.
We can then iterate this procedure (decomposing 
in the second step the function $|h|=\tilde g_2+\tilde h_2$ etc.) and 
obtain the above decomposition so that statements (i), (ii), (iii) 
hold. 

Also observe that
if $\La_n[|h|]\le 2^{-n}$  then $\fS_n[h]$ is contained 
in a dyadic cube of sidelength $2^{-n}$ and we thus know that
$\Theta_n[|h|]\La_n[|h|]=\int |h(x)|dx$.  This implies statement (iv).\qed
\enddemo

We now describe how the quantities of length and thickness are used in certain convolution 
estimates involving the measure $\mu$ and appropriate localizations $\mu^n$.
To define the localization  we choose 
a $C^\infty$ function $\phi$ with compact support in 
$\{x:|x|\le 1/2\}$ such that
$\int\phi(x) dx =1$ and such that 
$$\int\phi(x) (P(x)-P(0)) dx =0$$
for all polynomials of degree $\le d$.
Set $\phi_n(x)=2^{nd} \phi(2^n x)$ and let 
$$\mu^n=\phi_n*\mu. \tag 2.22$$

\proclaim{Lemma 2.4} Let $f$ be supported on a set of diameter at most 
$10$.
Then $$\meas(\supp(\mu^n *f))\lc \Lambda_n[f].$$
\endproclaim

\demo{\bf Proof}
Note that if $Q$ is a cube with center $x_Q$ and  sidelength $l(Q)$ 
 with $2^{-n}\le l(Q)\le 100$ 
and $f_Q$ is supported in $Q$ then
$\mu^n*f_Q$ is supported 
on  the $x_Q$-translate of a  tubular neighborhood
of $\Sigma$
of width $O(l(Q))$, thus on a set of  measure $O(l(Q))$. The assertion follows by
working with an efficient cover of the support of $f$ arising
 from the definition of $\Lambda_n$.\qed
\enddemo

The quantity $\Theta_n[f]$ can be used to estimate the $L^2$ norm of 
the support
$\mu^n *f$ provided that one has a lower bound for the curvature.
To make this precise we
 first prove a slight  variant of  an observation in \cite{7}.

\proclaim{Lemma 2.5} Let  $\psi$ be a real valued 
$C^\infty$ function on $[-1,1]^d$, so that
$\sup_{|\alpha|\le 3}|\partial^\alpha\psi(x)|\le A_3$; here $A_3\le 1$. 
Suppose  $|\det \psi''(y_0)|\ge \beta$ and  
$Q\subset [-1,1]^{d-1}$ is a $d-1$ dimensional cube
 of sidelength
$\eps_1 \beta$, containing $y_0$, here $\eps_1\le [10(d-1)^4 A_3]^{-1}$.

Let $\chi$ be a $C^\infty$ function supported on $Q$ so that the 
inequalities
$\|\partial^\alpha \chi\|_\infty\le 
c_\alpha (\eps_1\beta)^{-|\alpha|}$ hold.
Define the measure  $\nu$ by
$$\inn{\nu}{f}= \int\chi(y') f(y', \psi(y')) dy'$$
and define the reflection
$\inn{\widetilde {\nu}}{f}=\inn {\nu}{f(-\cdot)}$.

Then there are constants $C_\alpha$ so that
$$|\partial_x^\alpha[\nu*\widetilde {\nu}](x)|\le
C_\alpha\beta^{d-3-2|\alpha|}|x|^{-1-|\alpha|}.
$$
\endproclaim

\demo{\bf Proof}
We assume that $d\ge 3$ but after  notational modification the 
proof applies also to the case $d=2$.
Since $\nu*\tilde \nu$ does not change if we translate the measure
we may assume that $y_0=0$.

We compute 
$$\align 
\inn{\nu*\widetilde {\nu}}{f}&=\iint f(x-y) d\nu(x) d\nu(y)
\\&= \sum_k
\int\chi(u'+y')\chi(u') \zeta_k(u')
f(u', \psi(y'+u')-\psi(y')) dy' du':= \sum_k I_k(f)
\endalign
$$
where the $\zeta_k$ form a partition of unity  
on the unit sphere in $\Bbb R^{d-1}$ which is extended  to a homogeneous function of degree $0$. We assume that the restriction of 
$\zeta_k$  to the unit sphere 
is supported on a set of diameter  $\le \eps_1 \beta$ and the 
summation is over $O((\eps_1\beta)^{1-d})$ terms.
The $\zeta_k$ satisfy the natural estimates
$$|\partial^{\alpha} \zeta_k(u')|\le C_\alpha 
(\eps_1\beta)^{-|\alpha|}|u'|^{-|\alpha|}.$$
Note that in the integral defining $I_k$ the variables $u'$ and 
$y'$ are restricted to a ball of radius $\lc \eps_1\beta$ 
and $u'$ is further restricted to a sector with solid angle
$\eps_1\beta$.

Now note  that by $|\partial^2_{x_ix_j} \psi|\le A_3$, 
 $|\det \psi''(0)|\ge \beta$ and Cramer's rule  we have 
$$|u'|\le \beta^{-1}(d-1)^2 A_3^{d-1} |\psi''(0) u'|. \tag 2.23$$
We now pick a unit vector $\theta_k\in \supp \zeta_k$.

Let $$\fv_{k}=\frac{\psi''(0)\theta_k}{|\psi''(0)\theta_k|}$$ and 
let $\fv_{k,2},\dots ,\fv_{k,d-1}$ 
be an orthonormal basis of the orthogonal complement of $\bbR \fv_{k}$,
and with $t''=(t_2,\dots, t_{d-1})$ define
$\fw_k(t'')=\sum_{i=2}^{d-1}t_i\fv_{k,i}$.
Now write $y'=\fw_k(t'')+t_1 \fv_{k}$ and
we get
$$
I_k(f)=
\int_{t''}\int_{u'}
\int_{t_1} \chi(u')\zeta_k(u')
\chi(u'+ \fw_k(t'')+t_1 \fv_k)
f(u', \Psi_k(t_1, t'',u')) dt_1du' dt''
$$
where
$$
\align
\Psi_k(t_1, t'',u'))&=
 \psi(\fw_k(t'')+t_1 \fv_k+u')-
 \psi(\fw_k(t'')+t_1 \fv_k)
\\&=\biginn {u'}{\int_0^1 \nabla 
 \psi(\fw_k(t'')+t_1 \fv_k+su') ds}.
\endalign
$$
We  wish to change variables in the inner $t_1$-integral.
Observe that
$$\align
\frac d{dt}
\Psi_k(t_1, t'',u')=&
|u'| \inn{\theta_k}{\psi''(0)\fv_k}
\\&
+|u'|\int_0^1
\biginn{\theta_k}{ \big[\psi''(\fw_k(t'')+t_1 \fv_k+su')
- \psi''(0)\big] \fv_k} ds
\\
&+|u'|\int_0^1
\biginn{\frac{u'}{|u'|}-\theta_k}{\psi''(
\fw_k(t'')+t_1 \fv_k+su')  \fv_k } ds
\\=&
|u'||\psi''(0)\theta_k| +e_1(t_1,t'',u')+ e_2(t_1,t'',u')\tag 2.24
\endalign
$$ where by our assumption on the third derivatives
the error term $e_1$  is bounded by $ 2(d-1)^2A_3 
\eps_1 \beta |u'| $, and since $u'\in \supp \zeta_k$ the error term
$e_2$ 
is bounded by $(d-1)^2A_3 
\eps_1 \beta |u'| $.
The main term is 
$|u'||\psi''(0)\theta_k| \ge |u'|\beta (d-1)^{-2} A_3^{1-d}$ and thus the  derivative $\partial_t\Psi_k$   is 
single signed  and of size   $\approx \beta|u'|$. Therefore we may
perform the change of variables  
$t_1\mapsto u_d=\Psi_k(t_1, t'', u')$
with inverse $t_{1}^{k}(u_d; u',t'')$ and obtain
$$\align 
\inn{\nu*\widetilde {\nu}}{f}&=\sum_k
\iiint
f(u', u_d) H_k(u',u_d,t'') du_d du' dt''  \tag 2.25\endalign
$$
where
$$H_k(u',u_d,t'')=\frac{\zeta_k(u')
\chi(u')\chi(u'+\fw_k(t'')+t_1\fv_k)}
{|\partial_t \Psi_k(t_{1}^{k}(u_d;u',t''),t'',u')|
}.
$$
We have the estimate   $$|H_k(u',s, t'')|\lc\beta^{-1}|u'|^{-1}$$ and
  $H_k(u',u_d, t'')$ vanishes if 
$|u'|\ge C|u_d|$ or $|u'/|u'|-\theta_k|\ge \eps_1\beta$ or
  $|t''|\ge \beta$. Integrating in $t''$ yields a factor of 
$O(\beta^{d-2})$
and since $\sum_k \zeta_k (u')=O(1)$ 
we obtain the claimed estimate for  $\alpha=0$.  The estimates for the derivatives follow by a straightforward  examination of the 
derivatives of $t_1^k(u_d;u',t'')$ and 
applications of the chain rule. We omit the details.\qed
\enddemo

Now let $\phi_n$ be as in (2.22).
\proclaim{Lemma 2.6}
There is a small constant $\eps_1$ depending only on $\Sigma$  so that the following holds for
 $\beta\le 1$. 

Let $\chi\in C^\infty_0$ is supported on a set of diameter $\eps_1 \beta$ and suppose that
 the support of $\chi$ contains a point $P$ on $\Sigma$ where the Gaussian curvature 
satisfies $|K(P)|\ge \beta$.
Let
$\nu^n=\phi_n*\mu$.
Suppose that $f$ is supported on a set of diameter $1$. Then
$$
\|\widetilde{\nu}^n*\nu^n*f\|_\infty\lc \beta^{d-3} 
(1+n)\Theta_n[f]. 
$$
\endproclaim
\demo{\bf Proof} After localization and a change of variable we may 
reduce to the situation of Lemma 2.5.

Notice that $|\nu^n(x)|\lc 2^{n}$ since $\nu$ is a density on a hypersurface.
 By Lemma 2.5  we have
$$
|\widetilde {\nu}^n*\nu^n*f(x)|\lc \beta^{d-3}
 \int \min\{2^{n}, \frac 1{|x-y|}\}|f(y)|dy
$$
and we observe that
$$
 \int_{ |x-y|\le 2^{-n}} 2^{n}|f(y)|dy \le 2^d \Theta_n[f]
$$
and
$$ \int_{2^{-\ell}\le |x-y|\le 2^{-\ell+1}} \frac 1{|x-y|}|f(y)|dy \le
2^{d+1} \Theta_n[f], 
\quad 0\le \ell\le n.
$$
The asserted estimate follows by summing over $\ell=0,\dots,n$.\qed
\enddemo

Finally we  also  need the behavior of the quantities of length and thickness under nonisotropic dilations.
Here we will have to compare isotropic dilations to nonisotropic ones.
Let $\tau=\text{\sl trace}(P)$ and denote by $\la_j$ the eigenvalues of $P$. 
Then we may choose positive constants  $a$, $A$  so that
$$a<\text{Re}(\la_j)<A<\tau.\tag 2.26
 $$
Then there are positive constants $c_1\le C_1$ so that
 for all $x$ 
$$c_1t^a|x|\le |t^P x|\le C_1 t^A|x|, \qquad t\ge 1.
\tag 2.27
$$

\proclaim{Lemma 2.7} Suppose that  $f$ is integrable and
 vanishes in the complement of a compact set.

Then there is a constant $C$ depending only on the dilation group 
and the dimension, so that 

$$\Th_n[f(\delta_j \cdot)]\le C 2^{-j(\tau-A)}\Th_n[f],\quad\text{ if } j\ge 0\tag 2.28$$  
and
$$\La_n[f(\delta_{-m} \cdot)]\le C 2^{Am}\La_n[f], \quad\text{ if } m\ge 0. \tag 2.29$$

\endproclaim

\demo{\bf Proof}
Let $j\ge 0$ and let $Q$ be a dyadic cube of sidelength $l(Q)\ge 2^{-n}$. 
Then $\delta_{j} Q$ is contained in the union of at most $2^d$ 
dyadic cubes $\{q_i\}$, of sidelength
$\approx 2^{j A} l(Q)$. Thus
$$
\align
&l(Q)^{-1}\int_Q|f(\delta_j x)| dx=2^{-j\tau} l(Q)^{-1}\int_{\delta_j Q}|f(u)|du
\\
&\le 2^{-j\tau}\sum_{i} C'(2^{-Aj}l(q_i))^{-1}
\int_{q_i}|f(u)|du
\le C' 2^d 2^{-j(\tau-A)}\Theta_n[f].
\endalign
$$
If we take the supremum over all dyadic cubes we obtain  (2.28).

Next let $m\ge  0$. Let $Q_1,\dots, Q_N$ be a cover of $\fS_n(|f|)$. Let $Q_i^*$ be 
the double cube  (dilated with respect to the  center of $Q_i$).

Now $\fS_n(|f|)=\cup_{\nu=1}^{M_1} R_\nu$ where the $R_\nu$ are dyadic $2^{-n}$ cubes 
with center $x_\nu$ 
on which the expectation $E_n[|f|]$ does not vanish.
Let $R^*_{\nu,m}$ be the union of dyadic cubes of sidelength $2^{-n}$ 
 which intersect 
$\delta_m R_\nu$.  Then 
$\fS_n(|f(\delta_{-m}\cdot)|)$ 
is contained in $\cup_{\nu=1}^{M_1} R^*_{\nu,m}$.

Since  $m\ge 0$ each $R^*_{\nu,m}$ is contained in a $2$-dilate 
of $\delta_m R_\nu$ relative to the center $\delta_m x_\nu$. Thus  the union of the
$ R^*_{\nu,m}$ is contained in the union of the
 dilates $\delta_m Q_i^*$. Each 
$\delta_{m}Q_i^*$ is contained in  no more than
$4^d$ dyadic cubes of sidelength $ 2^{[m A+3]}l(Q_i)$.
Consequently
$$\La_n[f(\delta_{-m} \cdot)]\le 
C 2^{A m}\sum_{i=1}^N l(Q_i).
$$
If we work with an efficient cover of $\fS_n(|f|)$ we obtain (2.29).\qed
\enddemo

\head{\bf 3. Preliminary Calder\'on-Zygmund reductions }\endhead

We shall begin with some reductions from 
standard Calder\'on-Zygmund theory.
The estimates  in this section together with a trivial $L^1$ estimate
will only  imply the known  weak-type $L\log L$ inequality
(see Corollary 3.1 below) but they apply to more general operators 
than those discussed in 
the introduction.

In this section we shall  assume that the measure 
$\mu$  satisfies
$$|\widehat {\mu}(\xi)|
\lc (1+|\xi|)^{-\gamma}
\tag 3.1
$$
for some positive $\gamma$ (without loss of generality 
$\gamma\le (d-1)/2$).

When estimating the singular integral operator (1.5) we shall 
 assume the  additional cancellation condition  (1.4).
We note that  the original hypothesis 
of the curvature not vanishing to infinite 
order implies an estimate  (3.1) for some $\gamma>0$, by 
an application of van der 
Corput's lemma.

We shall apply a nonisotropic version 
of Calder\'on-Zygmund theory (see \cite{10}, \cite{16}).
Let $\rho$ be a homogeneous distance function which satisfies 
 $\rho(t^Px)=t\rho(x)$ for all $x$ 
 and $\rho(x)=1$ if $|x|=1$. 
If $x_0\in \Bbb R^d$ and $\rho_0>0$ then we set
$B(x_0,\rho_0)=\{x:\rho(x-x_0)\le \rho_0\}$ and we refer to $B(x_0,\rho_0)$ as the ball with center $x_0$ and $\rho_0$ (see \cite{17} for a discussion of such distance functions). Notice that
$$B(x_0,\rho_0)=\{x:|\rho_0^{-P}(x-x_0)|\le 1\}.$$
We note that 
$|x|^{1/a}\lc \rho(x)\lc |x|^{1/A}$ if $|x|\le 1$ and
$|x|^{1/A}\lc \rho(x)\lc |x|^{1/a}$ if $|x|\ge 1$, see (2.26/27) above.

Let $M_{HL}$ be the analogue of 
the Hardy-Littlewood maximal function associated to the family of these 
nonisotropic balls, {\it i.e.} 
$M_{HL}f(x)=\sup_{x\in B}|B|^{-1}\int_B|f(y)| dy$ where
 the supremum  is taken over all balls 
$B=B(x_0,\rho_0)$ which contain $x$.

We now fix $\alpha>0$ and define
$\Omega=\{x: M_{HL} f> \alpha\}$
and thus
$$|\Omega|\lc \alpha^{-1} \|f\|_1.
$$
By an analogue of the Lebesgue differentiation theorem we also know that
 $|f(x)|\le \alpha$ for all $x\in \bbR^d\setminus\Omega$.

The Calder\'on-Zygmund decomposition is based on a Whitney type 
decomposition. According to   \cite{16, p.15} there 
are constants $K_1>1$, $K_2>2$, $K_3>1$
(depending only on the distance function $\rho$),
 and a sequence of balls $B_1,\dots, B_j,\dots$, with $B_j=B(x_j,\rho_j)$, and a sequence $\fW$  of measurable sets (`generalized Whitney cubes')
$w_1,\dots, w_j,\dots$, 
 so that the following properties are satisfied:

\roster

\item"{{\it (a)}}"
The $B_j$ are pairwise disjoint.
\item"{{\it (b)}}"
If $B_j^*=B(x_j,K_1\rho_j)$ then the numbers $K_1\rho_j$  belong to 
$\{2^j:j\in \bbZ\}$ 
and  $\bigcup_j B_j^*=\Omega$. Moreover each $x\in \Omega$ is
contained in no more than $K_3$ of the balls $B_j^*$.

\item"{{\it (c)}}" $B_j\subset w_j\subset B_j^*$ 

\item"{{\it (d)}}"  The $w_j$ are pairwise disjoint, and we have 
$\bigcup w_j=\Omega$.

\item"{{\it (e)}}"
If $B_j^{**}=B(x_j,K_2\rho_j)$ then 
$ B_j^{**}\cap (\Bbb R^d\setminus\Omega)\neq\emptyset$.

\item"{{\it (f)}}" Each $B_j^{**}$ is contained in
$\Omega^*=\{x: M_{HL}(\chi_\Omega)> (10 K_2)^{-\tau}\}$ and thus
$$\meas(\Omega^*)\lc  \alpha^{-1}\|f\|_1\lc  \int \Phi(|f|/\alpha) dx.\tag 3.2$$
\endroster

We thus get a 
decomposition 
$ f=g+\sum_{w\in \fW} f_w$ where $f_w(x)=f(x)$  if $x\in w$ and
 $|f(x)|>\alpha$ and
$f_w(x)=0$ otherwise; moreover  $|g(x)|\lc  \alpha$ and
$|w|^{-1}\int |f_w|dx\lc \alpha$ for each $w$.
The sets $w$ play the role of the usual Whitney cubes.
For each $w\in \fW$ we assign a point $x_w$ and an integer $r(w)$ 
by setting
$x_{w_j}=x_j$ and $r(w_j)=\log_2(K_1\rho_j)$.

In what follows we choose 
$c>0$ small, specifically the choice
$$c< \frac12 \min\{1,\gamma\} \tag 3.3 $$ works.
We then further decompose $f_w$ by setting
$$f^n_w(x)= f_w(x)\quad\text{  if }\quad 2^{c (n-1)}\alpha < |f_w(x)|\le2^{c n}\alpha.$$
Observe that $f_w=\sum_{n=1}^\infty f^n_w$ and 
$$\sum_{n=1}^\infty \frac {1}{|w|}\int |f^n_w(x)| dx \lc \alpha.$$
We also let 
$$\align
g^n_w(x)&=  \chi_w(x) \frac{1}{|w|}\int_w f^n_w(y) dy, 
\\
b^n_w(x)&= f^n_w(x)- g^n_w(x),
\endalign
$$
and
$$
g^n(x)=\sum_w 
g^n_w(x), \qquad 
b^n(x)= \sum_w  b^n_w(x).
$$
Now
$$
\sum_{n=1}^\infty |g^n_w(x)|\le 
\frac{1}{|w|}\int_w \sum_{n=1}^\infty | f^n_w(y)| dy \, \chi_w(x)
\le  \frac{1}{|w|}\int_w | f_w(y)| dy \, \chi_w(x) \lc\alpha;
\tag 3.4
$$
moreover
$$
\sum_{n=1}^\infty |g^n(x)|\lc \alpha
\tag 3.5
$$
and 
$$
\sum_{n=1}^\infty \big[\|g^n_w\|_1+\|b^n_w\|_1\big] \lc \int_w |f(x)| dx \lc \alpha |w|.
\tag 3.6
$$

It will also be necessary to decompose the measure $\mu$  further. 
Let $\mu^n$ be the regularization defined in (2.22) and let
$$\mu^n_k(x)=2^{-k\tau }\mu^n(2^{-kP}x).$$

For our basic decomposition of the singular Radon transform 
we set $f^n=\sum_w f^n_w$ and using
$f=g+\sum_n f^n= 
g+\sum_n g^n +\sum_n b^n$ we split
$$
\sum_{k\in \bbZ} \mu_k *f= H_{I,1}+H_{I,2}+H_{I,3}+H_b
$$
where
$$
\aligned
H_{I,1}&= \sum_{k\in \bbZ} \mu_k* g
\\
H_{I,2}&=\sum_{k\in \bbZ} 
\sum_{n\ge 1}
(\mu_k-\mu_k^n)*  f^n
\\
H_{I,3}&=\sum_{k\in \bbZ} \sum_{n\ge 1} \mu_k^n *g^n
\\
H_{b}&= \sum_{k\in \bbZ} \sum_{n\ge 1} \mu_k^n *b^n.
\endaligned
\tag 3.7
$$

A further decomposition is necessary for $H_{b}$. 
For given $n\ge 1 $, $l\in \Bbb Z$ we define
$$\aligned
 I^n_l&=[ln, (l+1)n)\\
 (I^n_l)^*&=[(l-1)n, (l+1+\frac{2}{a})n]\endaligned
\tag 3.8
$$
 and set 
$$B^n_l=\sum_{w: r(w)\in I^n_l} b^n_w.$$
We split $H_b=H_{II}+H_{III}$
where
$$
\aligned
H_{II}&= \sum_{n\ge 1}\sum_{l\in \bbZ} \sum_{k\in \bbZ\setminus (I^n_l)^*} \mu_k^n *B^n_l
\\
H_{III}&= \sum_{n\ge 1}\sum_{l\in \bbZ} \sum_{k\in (I^n_l)^*} \mu_k^n *B^n_l.
\endaligned
\tag 3.9$$
Note that $H_{II}$ is 
 the portion of $H_b$ where the scaling of the measures 
$\mu_k^n$ is 
very different from the scaling of the balls $w$, which enables us to use 
standard $L^1$ arguments in the complement of the set $\Omega^*$.
The difficult term to estimate is $H_{III}$.

We shall show  that
$$
\align
\sum_{i=1}^3\|H_{I,i}\|_2
&\lc \alpha^{1/2} \|f\|_1^{1/2} \tag 3.10.1
\\
\|H_{II}\|_{L^1(\bbR^d\setminus \Omega^*)}&\lc \|f\|_1
\tag 3.10.2
\endalign
$$
From (3.10.1/2) we get by Chebyshev's inequality
$$\align
\meas\Big(\Big\{x: \sum_{i=1}^3|H_{I,i}(x)|>\alpha/{10}\}\Big)
&\lc \alpha^{-2}\Big\|
 \sum_{i=1}^3|H_{I,i}|\Big\|_2^2
\\
&\lc \alpha^{-2}\Big[ \sum_{i=1}^3\|H_{I,i}\|_2\Big]^2 \lc \alpha^{-1} \|f\|_1
\tag 3.11
\endalign
$$
and
$$
\meas\big(\big\{x\in \bbR^d\setminus \Omega^*: |H_{II}(x)|>\alpha/10\}\big)\lc \alpha^{-1}\|f\|_1.
\tag 3.12$$

We now prove the $L^2$ bounds (3.10.1) using standard arguments.
The cancellation of $\mu=\mu^0$ implies that $\widehat {\mu^0}(\xi)=O(|\xi|)$ and since $\mu_0$ is smooth we get
$$
|\widehat {\mu^0}(\xi)|\lc \min\{|\xi|, |\xi|^{-N}\}
\tag 3.13
$$
for large  $N$.

Even without such  a  cancellation assumption the difference 
$\mu^n-\mu^{n-1}$
does have cancellation and using the decay assumption (3.1) on the 
Fourier transform of $\mu$ it 
is straightforward to check that for $m\ge 1$
$$|\widehat {\mu^m}(\xi)-\widehat {\mu^{m+1}}(\xi)|
\lc 2^{-m\gamma} \min\{2^{-m}|\xi|, 
(2^{-m}|\xi|)^{-N}\}.
\tag 3.14
$$
Indeed the left hand side of (3.14) is $\lc (1+|\xi|)^{-\gamma}
|\widehat \phi(2^{-m}\xi)-\widehat \phi(2^{-m-1}\xi)|$
and since $\widehat\phi(\eta)=1+O(|\eta|^d)$  we obtain the bound
$2^{-m\gamma}(2^{-m}|\xi|)^{d-\gamma}$ which yields the claim for 
$|\xi|\le 2^{m+1}$ since also 
 $d-\gamma>1$. For $|\xi|\ge 2^{m+1}$ we use 
that $|\widehat \mu^m(\xi)|\le C_N |\xi|^{-\gamma}(1+2^{-m}|\xi|)^{-N}$.

Since $\widehat {\mu^n_k}(\xi)=\widehat{\mu^n}(\delta_{k}^*\xi)$ we obtain using (3.13), (3.14)
that
$$\aligned
&\sum_{k\in \Bbb Z}|\widehat {\mu^0_k }(\xi)| 
\lc 1
\\&\sum_{k\in \Bbb Z}
|\widehat {\mu^m_k }(\xi)-\widehat{\mu^{m-1}_k}(\xi)| 
\lc 2^{-m\gamma}.
\endaligned
\tag 3.15
$$

We recover the well-known result that $T$ is $L^2$ bounded, and as a consequence of the last displayed inequality we also get
$$
\Big\|\sum_{k\in \bbZ} (\mu_k-\mu_k^n)*f
\Big\|_2
\lc \sum_{m=n}^\infty
\Big\|\sum_{k\in \bbZ} (\mu^{m+1}_k-\mu_k^m)*f\Big\|_2
\lc 2^{-n\gamma}\|f\|_2.
$$
Now clearly
$$
\big\|H_{I,1}\big\|_2^2=\Big\| \sum_{k\in \bbZ} \mu_k* g\Big\|_2^2
\lc \|g\|_2^2\lc \alpha \|f\|_1
$$
and using (3.13) and (3.14) we also obtain
$$\align
\big\|H_{I,2}\big\|_2^2&\le
\Big(\sum_{n\ge 1}\Big\|\sum_{k\in \bbZ}
(\mu_k-\mu_k^n)*  f^n\Big\|_2\Big)^2\lc 
\Big(\sum_{n\ge 1}2^{-n\gamma}\| f^n\|_2\Big)^2
\\&\lc \sum_{n\ge 1}2^{-n\gamma}\| f^n\|_2^2
\lc 
\sum_{n\ge 1}2^{-n\gamma} \| f^n\|_1 
 2^{c(n+1)}\alpha\lc \alpha \|f\|_1
\endalign
$$
by our choice of $c$ in (3.3).
Moreover
$$
\align
\big\|
H_{I,3}\big\|_2^2&=\Big\|\sum_{k\in \bbZ} \sum_{n\ge 1} 
\big(\mu_k^0+\sum_{m=0}^{n-1}(\mu_k^{m+1}-\mu_k^m) *g^n\big)\Big\|_2^2
\\&\le \Big(
\Big\|\sum_{k\in \bbZ} \mu_k^0 * \sum_{n\ge 1} g^n\Big\|_2
+\sum_{m=0}^\infty\Big \|\sum_{k\in \bbZ} 
(\mu_k^{m+1}-\mu_k^m) *\sum_{n>m} g^n\Big\|_2\Big)^2
\\&\lc\Big(
\sum_{m=0}^\infty 2^{-m\gamma} \Big\|
\sum_{n>m} g^n\Big\|_2\Big)^2
\lc \alpha \|f\|_1.
\endalign
$$

Finally we prove the $L^1$ bound (3.10.2). Suppose that $r(w)\in I^n_l$.
For  $k\ge \max (I_n^l)^*$ (thus $k-r(w)\ge 2n/a$) 
we  use the cancellation of $b^n_w$ and obtain with $y_w\in w$
$$\align 
\mu^n_k* b^n_w(x)&=
\int 2^{-k\tau} \big[\mu^n(\delta_{-k}(x-y))
-
\mu^n(\delta_{-k}(x-y_w))\big] b^n_w(y) dy
\\&=
2^{-k\tau}\int \inn{\delta_{-k}(y-y_w)}
{\nabla \mu^n(\delta_{-k}(x-y_w+s(y-y_w)))} b^n_w(y) dy
\endalign
$$
and since 
$|\delta_{-k}(y-y_w)|\lc
2^{(r(w)-k)a}$ for $y\in w$ and $\|\nabla \mu^n\|_1=O(2^n)$ we get
$$\int
|\mu^n_k* b^n_w(x)| dx 
\lc 2^n 2^{(r(w)-k)a}\|b^n_w\|_1 .
$$ 
Moreover notice that by our assumption that $\mu$ is supported in the unit ball we have
that $\mu^n_k* b^n_w$ is supported in $\Omega^*$ if $k<\min (I^n_l)^*$.

Thus
$$\align \|
H_{II}\|_{L^1(\bbR^d\setminus \Omega^*)}&\le \sum_{n\ge 1}\sum_{l\in \bbZ} \sum_{k\ge \max (I^n_l)^*} \|\mu_k^n *B^n_l\|_1
\\
&\lc \sum_{n\ge 1}\sum_{l\in \bbZ} \sum_{k\ge \max (I^n_l)^*} \sum_{r(w)\in I^n_l}
2^n 2^{(r(w)-k)a}\|b^n_w\|_1 
\\
&\lc \sum_{n\ge 1}\sum_{l\in \bbZ} 2^{-n} \sum_{r(w)\in I^n_l}\|b^n_w\|_1 
\lc \|f\|_1,
\endalign
$$
by the definition of $(I^n_l)^*$. Thus  (3.10.2) is proved.

A  decomposition  similar to (3.7), (3.9)  applies to the maximal operator
where no cancellation on $\mu$ is assumed. We have  
$$
\sup_k| \mu_k *f|\le  M_{I,1}+M_{I,2}+M_{I,3}+M_{II}+M_{III}
$$
where
$$
\aligned
M_{I,1}&= \sup_{k\in \bbZ} |\mu_k* g|
\\
M_{I,2}&=
\sum_{n\ge 1} \sup_{k\in \bbZ} 
|(\mu_k-\mu_k^n)*  f^n|
\\
M_{I,3}&=\sum_{n\ge 1} \sup_{k\in \bbZ} 
|\mu_k^n *g^n|
\\
M_{II}&= \sum_{n\ge 1}\sum_{l\in \bbZ} \sup_{k\in \bbZ\setminus (I^n_l)^*} |\mu_k^n *B^n_l|
\\
M_{III}&= \sum_{n\ge 1}\sum_{l\in \bbZ} \sup_{k\in (I^n_l)^*} |\mu_k^n *B^n_l|
\endaligned
\tag 3.16
$$

Concerning the $L^2$ boundedness we observe that
$\sup_k|\mu^0_k*f|$ is pointwise controlled by the Hardy-Littlewood
maximal function $M_{HL} f$, associated to the given dilation group.
Therefore 
$$\big \|\sup_k|\mu^0_k*f|\big\|_2\lc\|f\|_2.
\tag 3.17
$$
Again by Fourier transform arguments as above
$$\align
&\big \|\sup_k|(\mu^m_k -\mu^{m-1}_k)*f|\big\|_2
\lc
\Big \|\Big(\sum_k|(\mu^m_k -\mu^{m-1}_k)*f|^2\Big)^{1/2}\Big\|_2
\\
&\lc
2^{-m\gamma}\Big(\int\sum_{k}
|\widehat{\mu^m_k}(\xi) -\widehat{\mu^{m-1}_k}(\xi)|^2
|\widehat f(\xi)|^2 d\xi \Big)^{1/2}\lc
2^{-m\gamma}\|\widehat f\|_2\lc 2^{-m\gamma}\|f\|_2.
\endalign
$$
This shows that  we can repeat the arguments for $H_I$ above and get
$$
\sum_{i=1}^3\|M_{I,i}\|_2
\lc \alpha^{1/2} \|f\|_1^{1/2}. \tag 3.18
$$
In the definition of $M_{II}$ we may replace the $\sup$ over
${k\notin (I^n_l)^*}$ by the sum 
and the estimation 
is exactly the same as for $H_{II}$ above. This yields
$$\|M_{II}\|_{L^1(\bbR^d\setminus \Omega^*)}\lc \|f\|_1.
\tag 3.19
$$
We combine these estimates with (3.2)
 and we see that in order to prove Theorems 1.1 and 1.2 
we are left to prove
the inequalities
$$\align
&\meas\{x: |M_{III}|>\frac 45\alpha\} \lc \int\frac{|f(x)|}{\alpha}
\log\log(e^2+\frac{|f(x)|}{\alpha}) dx
\tag 3.20
\\
&\meas\{x: |H_{III}|>\frac 45\alpha\} \lc  \int\frac{|f(x)|}{\alpha}
\log\log(e^2+\frac{|f(x)|}{\alpha}) dx
\tag 3.21
\endalign
$$
This will be done in \S5 and \S6 below.

\remark{Weak type  $L\log L $ estimates}
We  note that   weak type $L\log L$ inequalities for $T$ and $\cM$ can be 
already obtained from
trivial $L^1$ estimates for $H_{III}$ and $M_{III}$.
Here  we are essentially reproving the result in \cite{4}.

\proclaim{Corollary 3.1}
Let $\mu$ be a compactly supported Borel measure satisfying
$$
|\widehat {\mu}(\xi)|\le C(1+|\xi|)^{-\gamma}.
$$
Then $M$ is of weak type $L\log L$.
If in addition the cancellation condition $\int d \mu(x)=0$ holds, then $T$ is of weak type $L\log L$.
\endproclaim

\demo{\bf Proof} Given our previous estimates we just have to estimate
the measure of the sets where $M_{III}>\alpha$ or $|H_{III}|>\alpha$.
We simply use Chebyshev's inequality and are left with 
estimating
$\alpha^{-1}\|M_{III}\|_1$ and 
$\alpha^{-1}\|H_{III}\|_1$, respectively.
Using that the $L^1$ norm of $\mu^n_k$ is uniformly bounded in $k,n$ we get
$$\align \|H_{III}\|_1&\le 
\sum_{n\ge 1}\sum_{l\in \bbZ} \sum_{k\in (I^n_l)^*} \|\mu_k^n *B^n_l\|_1
\lc
\sum_{n\ge 1}\sum_{l\in \bbZ} \sum_{k\in (I^n_l)^*}\sum_{r(w)\in I^n_l} \|b_w^n\|_1
\\&\lc
\sum_{n\ge 1} \sum_{l\in \bbZ}
 \sum_{r(w)\in I^n_l}n \|b_w^n\|_1 
\lc
\sum_{n\ge 1}  n\|f^n\|_1 \lc \int|f(x)|\log(e+\frac{|f(x)|}{\alpha}) dx
\endalign
$$
and the same argument applies to $M_{III}$.\qed
\enddemo
\endremark

\head{\bf 4. A stopping time argument}\endhead

In order to refine the previous estimates for
 $M_{III}$ and $H_{III}$ we need a further 
decomposition of $b^n_w$. Here we
use a stopping time argument  based on length $\Lambda_n$ (and
 thickness $\Theta_n$). The reader will note 
some  similarities with 
Christ's stopping time argument in \cite{2}.

In what follows $\fQ_0$ will denote the set of dyadic unit cubes 
of the form 
$(n_1,\dots, n_d)+[0,1)^d$, $n_i\in \bbZ$.

\proclaim{Proposition 4.1}
For every $n$ 
and every $w$ with $r(w)\in I^n_l$
 there is a  decomposition
$$
b^n_w= \sum_{\kappa\in (I^n_l)^*} f^{n,\ka}_w
\tag 4.1
$$
so that the following properties are satisfied.

(i)
$$
\sum_{\kappa\in (I^n_l)^*} |f^{n,\ka}_w|= |b^{n}_w|.
\tag 4.2
$$

(ii) For every $q\in \fQ_0$, $\kappa\in (I^n_l)^*$
$$\Lambda_n\big[\sum_{r(w)<\ka} f^{n,\ka}_w(\del_\ka\cdot)\chi_q\big]
\le \alpha^{-1} \sum_{r(w)<\ka}\int_q|\fnkw(\del_\ka y)|dy.
\tag 4.3
$$

(iii) For every $q\in \fQ_0$, and for every $\kappa\in (I^n_l)^*$ and 
$s\ge 1$ with $\kappa+s\in (I^n_l)^*$,
$$
\Th_n\big[\sum_{r(w)\le \ka}\fnkw(\del_{\ka+s}\cdot)\chi_q\big]\le
 16(n+1)\alpha.
\tag 4.4
$$
\endproclaim

\demo{\bf Proof}
This is proved by an inductive construction.

We shall give a decomposition of 
$$\cG^0=\sum_{w:r(w)\in I^n_l} b^n_w;
$$
since the $w$ are disjoint this  will  yield a decomposition 
of each  $b^n_w$.
 Set 
$\ka^\mx_{n,l}=\max (I^n_l)^*$ and 
$\ka_j=\ka^\mx_{n,l}-j$. We shall establish  the following

\proclaim{\bf Claim} For  $N=0,1,\dots$ we can decompose
$$\cG^0= \sum_{j=0}^N [H^j+S^j] + G^N$$
so that

(i)  $G^{j-1}=H^j+S^j+G^j$ if $j\ge 1$

(ii)
$G^j=\sum_{q\in \fQ_0}\sum_{\nu=1}^{L(j,Q)} G^{j,q}_\nu$,
where $G^{j,q}_\nu$ vanishes in the complement of $\del_{\ka_j}q$ and 
$$\Theta_n[G^{j,q}_\nu(\delta_{\ka_j}\cdot)]\le 8\alpha.$$

Moreover $$L(j,Q)\le n+1.$$

(iii) 
$$
\aligned
H^j(x)&=0 \qquad \text{if } x\notin \bigcup_{r(w)<\kappa_j}w
\\
S^j(x)&=0 \qquad \text{if } x\notin \bigcup_{r(w)=\kappa_j}w
\\
G^N(x)&=0  \qquad \text{if } x\notin \bigcup_{r(w)<\kappa_N}w.
\endaligned
$$

(iv) For each $q\in \fQ_0$,
$$\Lambda_n\big[H^j(\delta_{\ka_j}\cdot)\chi_q\big]\le \alpha^{-1}
\int_q|H^j(\delta_{\kappa_j}y)|dy.
$$

(v) For $\kappa>\kappa_j$, $\kappa\in (I^n_l)^*$ and each
$q\in \fQ_0$,
$$\Theta_n\big[H^j(\delta_\kappa\cdot)\chi_q\big]
+
\Theta_n\big[S^j(\delta_\kappa\cdot)\chi_q\big]
\le 16(n+1)  \alpha. 
$$

(vi) The functions  $G^j$, $G^{j,q}_\nu$, $H^j$, $S^j$ are nonnegative at $x$ (nonpositive) if and only if $f(x)$ is nonnegative (nonpositive).
\endproclaim

If we accept the claim
then 
in order to complete the proof of the proposition we observe that 
in the above statement $\kappa=\kappa_j=\kappa^\mx_{n,l} -j$ and thus 
we  merely have to define
$$
f^{n,\kappa}_w(x)=
\cases
H^{\kappa^\mx_{n,l} -\ka}
(x) \qquad&\text{ if } x\in w, r(w)<\kappa\le 
\kappa^\mx_{n,l}, 
\\
 S^{\kappa^\mx_{n,l} -\ka}(x) \qquad&\text{ if } x\in w,  r(w)=\kappa
\\
0&\text{ if }x\notin w \text{ or if } \kappa<r(w).
\endcases
$$
Then  (4.1) follows from (iii) and (4.2) from (4.1) and (vi).
(4.3) is a consequence of (iv) and (4.4) follows from (v).

\demo{Proof of the Claim} We argue by induction and assume that either 
$N=0$ or that  $N>0$ and  statements (i)-(vi) hold for all $j\le N-1$.

If $N=0$ we set $S^0=H^0=0$ and $G^0=\cG^0$.
If $N\ge 1$ we begin by defining functions $S^N$, $\cG^N$
where
$S^{N}(x)=G^{N-1}(x)$ if $x\in \bigcup_{r(w)=\kappa_N} w$ and $S^N(x)=0$ 
otherwise, and 
$\cG^N(x)= G^{N-1}(x)-S^N(x)$. Thus $\cG^N$ is supported on 
$\bigcup_{r(w)<\ka_N}w$ and coincides with $G^{N-1}$ there. 
Note that $\cG^N$ vanishes if 
$\kappa_j<\min I^n_l$ and the construction stops then.

We now use Proposition 2.3 to decompose for $q\in \fQ_0$
$$
\cG^N(\delta_{\ka_N}x)\chi_q(x)
 =\sum_{\nu=1}^L g^{N,q}_\nu+ h^{N,q}_L
$$
so that 
$\Theta_n[g^{N,q}_\nu]\La_n[h^{N,q}_{\nu-1}]\le 8\int|g^{N,q}_\nu|dx$
and
$h^{N,q}_L$ vanishes for  $L\ge n+1$.
Also the signs of the functions
$g^{N,q}_\nu$, $ h^{N,q}_L$ coincides with the sign of
$\cG^N(\delta_{\ka_N}(x))\chi_q(x)$ and we have 
$h^{N,q}_{\nu-1}= g^{N,q}_\nu+ h^{N,q}_\nu$ for $\nu\ge 1$ with
$\Lambda_n[h^{N,q}_\nu]\le
\Lambda_n[h^{N,q}_{\nu-1}]/2$.

Let 
$L(N,q)$ be the {\it minimal} integer $L$ so that
$$\Lambda_n[h^{N,q}_L]
\le \alpha^{-1}\int|h^{N,q}_L(y)| dy.\tag 4.5$$
Then $L(N,q)\le n+1$
(since 
$h^{N,q}_L$ vanishes for  $L\ge n+1$).

Now,
$\Lambda_n[h^{N,q}_{\nu-1}]\ge \alpha^{-1}\int|h^{N,q}_{\nu-1}(y)| dy$
for  $\nu\le L(N,q)$,
by the minimality of $L(N,q)$, 
 and since $|g^{N,q}_\nu|\le |h^{N,q}_{\nu-1}|$ we get
$$
\Th_n[g^{N,q}_\nu]\le 8
\frac{\int|g^{N,q}_\nu(y)|dy}{\Lambda_n[h^{N,q}_{\nu-1}]} \le
8\alpha.
\tag 4.6
$$

Now define
$G^{N,q}_\nu(x)= g^{N,q}_\nu(\delta_{-\ka_N}x)$, for $\nu\le L(N,q)$, and
$G^N(x)=\sum_{q\in \fQ_0} 
\sum_{\nu=1}^{L(N,q)}
G^{N,q}_\nu(x)$. Moreover
$H^{N,q}(x) =h^{N,q}_{L(N,q)}(\delta_{-\ka_N}x)$ and $H^{N}(x)=\sum_{q\in \fQ_0} H^{N,q}(x)$.
 Then the statement (vi) about  the sign of $G^{N,q}_\nu$, $G^N$ and $H^N$
holds. (iv) follows from (4.5).
Statements (i) and (iii) hold by construction, 
and the inequality for the thickness 
in (ii) holds by (4.6) by (4.6).

In view of (i), (vi) we also  have 
$|H^N|+|S^N|\le |G^{N-1}|\le |G^{N-s}|$  for $s\ge 1$ so  that
by statement (ii)  for $j\le N-1$ we get
$$
\align
&\Theta_n[H^N(\del_{\ka_N+s}\cdot)\chi_q]
+\Theta_n[S^N(\del_{\ka_N+s}\cdot)\chi_q]
= 
\Theta_n[H^N(\del_{\ka_{N-s}}\cdot)\chi_q]
+\Theta_n[S^N(\del_{\ka_{N-s}}\cdot)\chi_q]
\\&\le
2 \Theta_n[G^{N-s}(\delta_{\ka_{N-s}}\cdot)\chi_q]
\le 16 \sum_{\nu=1}^{L(N-s,q)}\Theta_n[G^{N-s,q}_\nu(\delta_{\kappa_{N-s}}\cdot)]
\le 16 L(N-s,q)\alpha\le 16(n+1)\alpha
\endalign
$$
This implies $(v)$ for $j=N$ and the Claim is proved.
\qed

\enddemo
\enddemo

\head{\bf 5. The main estimate for the maximal function}\endhead

We shall prove  the nontrivial estimate (3.20) for the 
maximal function,  assuming again  that the curvature assumption in 
the introduction is satisfied, and  prove the inequality
$$\meas\Big(\Big\{x:\sup_{k}\big|
\sum\Sb n,l\\k\in (I^n_l)^*\endSb  \mu_{k}^{n}*B^{n}_l
\big|
>\alpha\Big\}\Big) \le \int \Phi (|f|/\alpha) dx\tag 5.1
$$
with $\Phi(t)=t\log\log (e^2+t)$.

%

We use the decomposition in Proposition 4.1 and form an additional exceptional set $\cO_1$.
To define it we set  for $q\in \fQ_0$, $\kappa\in (I^n_l)^*$,
$$F^{n,l,\ka}_q(x)=\sum\Sb r(w)\in I^{n}_l\\ r(w)<\ka\endSb f^{n,\ka}_w(x) 
\chi_{q}(\delta_{-\ka}x).
\tag 5.2$$
and define
 
$$\cO_1=\bigcup_{n=1}^\infty \bigcup_{l\in \Bbb Z}
\bigcup\Sb \kappa\in (I^n_l)^* \endSb
\bigcup\Sb  q\in \fQ_0\endSb
\bigcup\Sb k\in (I^n_l)^*\\ k\le \ka\endSb 
\supp\big(\mu_{k}^n*F^{n,l,\ka}_{q}\big);
\tag 5.3
$$
moreover we define
$$\cO=\cO_1\cup \Omega^*\tag 5.4$$
where $\Omega^*$ is as in (3.2).

To estimate the measure of $\cO_1$ observe that
$\supp(\mu_{k}^{n}*F^{n,l,\ka}_{q})=\delta_k
\supp(\mu_{0}^{n}*[F^{n,l,\ka}_{q}(\delta_k\cdot)])$ and since 
for $k\le \ka$ the function 
$F^{n,l,\ka}_{q}(\delta_k\cdot)$ is supported in a set of bounded diameter
 we get by
 (2.29) and (4.3)
$$
\align
&\meas\big(
\supp(\mu_{k}^{n}*F^{n,l,\ka}_{q})\big)
= 2^{k\tau}
\meas\big(
\supp(\mu_{0}^{n}*[F^{n,l,\ka}_{q}(\delta_k\cdot)])\big)
\\&\lc 2^{k\tau} \Lambda_n[F^{n,l,\ka}_q(\delta_k\cdot)]\lc 2^{k\tau} 2^{(\ka-k)A}
\La_n [F^{n,l,\ka}_{q}(\delta_\ka\cdot)]
\\
&\lc 2^{k\tau} 2^{(\ka-k)A}  \alpha^{-1}\int|F^{n,l,\ka}_{q}(\delta_\ka y)|dy
\lc  2^{(k-\ka)(\tau-A)} \alpha^{-1}\int|F^{n,l,\ka}_{q}(y)|dy.
\endalign
$$

Thus,  we can sum a geometric series in $k\le \ka$ and obtain
$$\align \meas(\cO_1)&\lc
\sum_{n=1}^\infty  \sum_{l\in \Bbb Z}
\sum\Sb \kappa\in (I^n_l)^*\endSb
\sum_{q\in \fQ_0} \alpha^{-1}\int|F^{n,l,\ka}_{q}(y)|dy
\lc\sum_{n=1}^\infty  \sum_w
 \alpha^{-1}\int|b^{n}_{w}(y)|dy
\\&\lc\alpha^{-1}
\sum_{n=1}^\infty  \sum_w
 \int|f^{n}_{w}(y)|dy
\lc \alpha^{-1}\int|f(y)| dy
\tag 5.5
\endalign$$
and  the measure of $\cO=\cO_1\cup \Omega^*$ satisfies the same estimate.
Note that the contributions for $k\le \kappa$, $r(w)=\ka$ are also supported in  
$\cO$ since $\mu$ is assumed to be supported in the unit ball 
and thus
$$\bigcup_{n=1}^\infty \bigcup_{l\in \Bbb Z}
\bigcup\Sb w:r(w)\in I^n_l \endSb
\bigcup\Sb  k\le r(w)\endSb 
\supp\big(\mu_{k}^n*f^{n,r(w)}_{w}\big)\subset \Omega^*.
$$

It now remains to handle the contribution in the 
complement of $\cO$ which 
only involves the scales $k>\ka$ and contributions
for $r(w)\in I^n_l$ with $r(w)\le \ka$; to simplify the notation below we set
$$I^{n,\kappa}_l=\{r\in I^n_l: r\le \ka\}.$$

We shall first cut out a contribution from   'flat' parts of $\Sigma$.
We recall that the curvature does not vanish to infinite order on 
$\Sigma$ and therefore there is 
 a number  $\eta>0$ such that
$$
\int_{\Sigma} |K(x)|^{-\eta} d\sigma(x)<\infty.
\tag 5.6
$$
This is well known (for example, one may use an argument 
in \cite{16, p.343} to reduce to an inequality in one dimension where one can use
  H\"older's inequality and compactness).

By Chebyshev's inequality (5.6) implies that
$$
| \{x\in \Sigma: |K(x)|\le n^{-3/\eta}\}| \lc n^{-3}.
\tag 5.7
$$
Now we use a partition of unity to write 
$$\mu=\sum_{i\in \cJ^n}\nu^{i,n}$$
where each  $\nu^{i,n}$ is supported 
on a cube  $R_i$ 
of diameter $\eps_1 n^{-3/\eta}$ (here $\eps_1$ will be as in Lemma 2.6)
and the supports of the $\nu^{i,n}$ have bounded overlap, 
independent of $n$.
Note that then
$$\card (\cJ^n)\lc n^{3(d-1)/\eta}. 
\tag 5.8
$$

We  split the index set into disjoint subsets 
 as  $\cJ^n=\cJ_1^n\cup \cJ_2^n$ where
 $\cJ_2^n$ consists of 
all $i\in \cJ$ with the property that 
$|K(x')| \le n^{-3/\eta}$ for all $x'\in \supp R_i$. 

Then  by (5.7) we have that the sum of the  total variations of the 
$\nu^{i,n}$,  for which  $i\in \cJ^n_2$, 
satisfies the bound
$$
\sum_{i\in \cJ_2^n} \big\|
\nu^{i,n}
\big\|\lc n^{-3}.
$$
Let $$\mu^{i,n}=\nu^{i,n}*\phi_n$$ and
$\mu^{i,n}_k=2^{-k\tau}\mu^{i,n}(2^{-kP}\cdot)$.

Since the cardinality of $(I^n_l)^*$ is $O(n)$  and 
$\sum_{i\in \cJ_2^n} \|
\mu^{i,n}_k\|_1=O(n^{-3})$
 the contribution of the measures 
$\sum_{i\in \cJ_2^n} \mu^{i,n}_k$, $k\in (I^n_l)^*$ can be handled by a 
straightforward $L^1$ estimate:
$$\align
&\meas\Big(\Big\{x:\sup_{k}\Big|
\sum\Sb n,l\\k\in (I^n_l)^*\endSb
\sum\Sb\ka\in (I^n_l)^*\\ \ka\le k\endSb  \sum_{i\in \cJ_2^n}
|\mu^{i,n}_k*\sum\Sb r(w)\in I^{n,\ka}_l\endSb f^{n,\ka}_w\Big|>
\alpha/10\Big\}\Big)
\\&\lc \alpha^{-1}
\Big\|
\sum\Sb n,l\endSb\sum \Sb k\in (I^n_l)^*\endSb
\sum\Sb\ka\in (I^n_l)^*\\\ka\le k\endSb  \sum_{i\in \cJ_2^n}
|\mu^{i,n}_k*\sum\Sb r(w)\in I^{n,\ka}_l\endSb f^{n,\ka}_w|\Big\|_1
\\&\lc \alpha^{-1}
\sum\Sb n,l\endSb
\sum\Sb\ka\in (I^n_l)^*\endSb
\sum \Sb k\in (I^n_l)^*\endSb \sum_{i\in \cJ_2^n}\big\| 
\mu^{i,n}_k\big\|_1\sum\Sb r(w)\in I^{n,\ka}_l\endSb \|f^{n,\ka}_w\|_1
\\
&\lc \alpha^{-1}
\sum\Sb n,l\endSb
\sum\Sb\ka\in (I^n_l)^*\endSb
n^{-2}\sum\Sb r(w)\in I^{n,\ka}_l\endSb \|f^{n,\ka}_w\|_1
\lc\alpha^{-1}\|f\|_1.
\tag 5.9
\endalign
$$

Next choose a large constant $C_0$; specifically  the choice  $$C_0\ge 
\frac{100}{a}(1+\frac{d}{\eta})\max\{1,\frac{A}{\tau-A}\} +10+ \log_2 \big(\frac{C_1}{c_1}\big)
\tag 5.10$$ 
will work where $c_1\le  C_1$ are as in (2.27).
Then the  contribution for the scales $\ka\le k\le \ka + C_0 \log n$ is 
also handled by an $L^1$ estimate:

$$\align
&\meas\Big(\Big\{x:\sup_{k}\Big|
\sum\Sb n,l\\ k\in (I^n_l)^*\endSb
\sum\Sb\ka\in (I^n_l)^*\\ \ka\le  k\le \ka +C_0 \log n\endSb  \sum_{i\in \cJ_1^n}
\mu^{i,n}_k*\sum\Sb r(w)\in I^{n,\ka}_l\endSb 
f^{n,\ka}_w(x)\Big|>\alpha/10\Big\}\Big)
\\&\lc \alpha^{-1}
\Big\|
\sum\Sb n,l\endSb\sum \Sb k\in (I^n_l)^*\endSb
\sum\Sb\ka\in (I^n_l)^*\\ \ka\le k\le \ka+C_0\log n\endSb  \sum_{i\in \cJ_1^n}
|\mu^{i,n}_k*\sum\Sb r(w)\in I^{n,\ka}_l\endSb f^{n,\ka}_w|\Big\|_1
\\&\lc \alpha^{-1}
\sum\Sb n,l\endSb \sum\Sb\ka\in (I^n_l)^*\endSb
\sum_{w:r(w)\in I^{n,\ka}_l}
\sum \Sb \ka\le k\le \ka+C_0\log n\endSb \|f^{n,\ka}_w\|_1
\\
&\lc \alpha^{-1}
\sum\Sb n,l\endSb\sum\Sb\ka\in (I^n_l)^*\endSb\sum_{r(w)\in I^{n,\ka}_l}
\log n \|f^{n,\ka}_w\|_1
\\&\lc\alpha^{-1}
\sum\Sb n\endSb
\log n \|f^n\|_1
\lc \int\frac{|f(x)|}{\alpha} \log\log\big(e^2+
\frac{|f(x)|}{\alpha}\big) dx
\tag 5.11
\endalign
$$

It remains to show
$$\meas\Big(\Big\{x:\sup_{k}\Big|
\sum\Sb n,l\\ k\in (I^n_l)^*\endSb
\sum\Sb\ka\in (I^n_l)^*\\  k> \ka +C_0 \log n\endSb  \sum_{i\in \cJ_1^n}
\mu^{i,n}_k*\sum\Sb r(w)\in I^{n,\ka}_l
\endSb f^{n,\ka}_w(x)\Big|>\alpha/10\Big\}\Big) \lc \alpha^{-1}\|f\|_1
\tag 5.12
$$
and this will be accomplished by proving $L^2$ estimates.

\subheading {Reintroducing  cancellation} The decomposition in (4.1) was needed to exploit the geometry of the exceptional set; however we paid the price of destroying the cancellation properties of the $b^n_w$. As the 
information on the support of the $f^{n,\ka}_{w}$ has been used and is not needed anymore for the scales $k>\ka+C_0 \log n$ we shall now
modify the functions $f^{n,\ka}_w$  to reintroduce some cancellation.
Namely let $\{P_i\}_{i=1}^{M_d}$ be an orthonormal basis 
of the space of polynomials  of degree
$ \le d$
 on the unit ball  $\{x:|x|\le 1\}$
and for given $w$ define
 the projection operator $\Pi_w$ by
$$\Pi_w[h](x)= \chi_w(x)\sum_{i=1}^{M_d}  P_i(\delta_{-r(w)}(x-x_w))
 \int_w h(y) P_i(\delta_{-r(w)}(y-x_w)) 2^{-r(w)\tau}\,dy.
$$
Note that
$$\big|\Pi_w[h](x)\big|\le C \frac{1}{|w|}\int_w|h(y)|dy
\tag 5.13
$$
where $C$ is independent of
$h$ and $w$.

Let 
$$\align
g^{n,\ka}_w(x) &=\Pi_w [f^{n,\ka}_w](x),
\\b^{n,\ka}_w(x) &= f^{n,\ka}_w(x)-g^{n,\ka}_w(x),
\endalign
$$
so that $b^{n,\ka}_w$ vanishes off  $w$ and for polynomials $p$
$$\int_w b^{n,\ka}_w(x) p(x) dx =0 
\quad \text{ if } \text{deg}(p)\le d.
\tag 5.14 $$

We observe that since the  $w$'s are generalized  Whitney cubes for $\Omega$ (see \S3), we have
$$
\sum_{n,\ka}\big|\Pi_w [f_w^{n,\ka}](x)\big|  \lc
\chi_w(x) \frac{1}{|w|}\int_w |f(x) | dx \lc \alpha ;
\tag 5.15$$
moreover by (5.13) 
$$\sum_{n,\ka}\big[\|b^{n,\ka}_w\|_1 +\|g^{n,\ka}_w\|_1 \big]
\lc 
\sum_{n,\ka}\|f^{n,\ka}_w\|_1 \lc
\int_w|f(x)| dx.\tag 5.16
$$

Now (5.12) will follow 
from
$$\align
&\Big\|
\sup_{k}\Big|
\sum\Sb n,l\\ k\in (I^n_l)^*\endSb
\sum\Sb\ka\in (I^n_l)^*\\  k> \ka +C_0 \log n\endSb  \sum_{i\in \cJ_1^n}
\mu^{i,n}_k*\sum\Sb r(w)\in I^{n,\ka}_l\endSb g^{n,\ka}_w\Big|\Big\|_2^2\lc \alpha\|f\|_1
\tag 5.17
\\&\Big\|
\sup_{k}\Big|
\sum\Sb n,l\\ k\in (I^n_l)^*\endSb
\sum\Sb\ka\in (I^n_l)^*\\  k> \ka +C_0 \log n\endSb  \sum_{i\in \cJ_1^n}
\mu^{i,n}_k*\sum\Sb r(w)\in I^{n,\ka}_l\endSb b^{n,\ka}_w\Big|\Big\|_2^2\lc \alpha\|f\|_1.
\tag 5.18\endalign
$$

The estimation  (5.17) is straightforward.
If $d\sigma$ denotes surface measure on $\Sigma$ and $d\sigma_k$ the dilate 
$2^{-k\tau}d\sigma(\delta_{-k}\cdot)$ then
the maximal function $$Mf(x)=\sup_{k\in \bbZ} |d\sigma_k* f|$$ defines a 
bounded operator on $L^2$. By the positivity of this maximal operator 
the left side of (5.17) is bounded by a constant times
$$
\Big\|M_{HL}M\big[\sum\Sb n,l\endSb
\sum\Sb\ka\in (I^n_l)^*\endSb 
\sum\Sb w\endSb |g^{n,\ka}_w|\big]\Big\|_2^2
\lc\alpha\sum\Sb n,l\endSb
\sum\Sb\ka\in (I^n_l)^* \endSb \sum\Sb w\endSb 
\big\|g^{n,\ka}_w\big\|_1 \lc \alpha \|f\|_1;
$$
here we used (5.15/16).

For the remainder of this section we prove (5.18). 

We first replace the $\sup$ in $k$ 
by an $\ell^2$ sum and then, for fixed $k$, we apply Schwarz' inequality 
in the form 
$[\sum_n |a_n|]^2\lc \sum |na_n|^2$.
Next we observe  that for fixed $n$ the 
number $k$ is contained  in at most $3+2/a$ of the intervals 
 $(I^n_l)^*$.
Then we apply Schwarz' inequality
 for the sim in  $\ka$  yielding a factor of $O(n)$ and 
for the sum in 
$i$ yielding a factor of  $O(n^{3(d-1)/\eta})$. 
Finally we group the sum over $w$ into groups for which
$r(w)=r$, $r\in I^n_l$ and apply 
Schwarz' inequality in $r$ which yields one more factor of $O(n)$. Thus we see that 
the left side of (5.18) is dominated by a constant times

$$
\sum\Sb k,n,l\\k\in (I^n_l)^* \endSb \sum\Sb \ka: \ka<\\
k-C_0\log n\endSb
\sum_{i\in \cJ_1^n}
\sum\Sb r\in I^{n,\ka}_l\endSb
n^{(4+\frac{3(d-1)}{\eta})}
\Big\|\mu^{i,n}_k*\sum\Sb r(w)=r\endSb b^{n,\ka}_w\Big\|_2^2
\tag 5.20
$$

We note that the some of the applications  of Schwarz' inequality
 above are not really  necessary but it turns out that the polynomial factors in $n$ are irrelevant in the range
$\ka<k-C_0\log n$.

Now, for fixed $\ka, k$,  define
$$M(\ka, k) = \big[k-(k-\ka)\frac{a}{2 A}+\log_2\frac{C_1}{c_1}+2 ]
\tag 5.21$$
where $[v]$ denotes the largest integer $\le v$.
Note that for $\ka<k-C_0\log n$ we have $M(\ka,k)<k$.
Let  $\fR(\ka,k)$ be the collection of dilates
$\delta_{M(\ka,k)}q$, where  $q\in \fQ_0$.
For each $w$ with $r(w)=r\le \ka$ we assign $R\in \fR(\ka, k)$ so that 
$w\cap R\neq\emptyset$.
We write $R=R_{\ka,k}(w)$ or simply $R=R(w)$ if the dependence on $k,\ka$ is clear.

 Let  $\widetilde \fR(\ka,k)$ be a subcollection of $\fR(\ka, k)$ 
with the property that
if $R, R'\in \fR(\ka, k)$, $R\neq R'$  and $R=\delta_{M(\ka,k)}q$,
$R'=\delta_{M(\ka,k)}q'$ then $\dist (q,q')>10$.

We shall show  for fixed $n$, $l$, $k\in (I^n_l)^*$, $\ka\in (I^n_l)^*$,
 $r\in I^{n,\ka}_l$ that
$$\Big\|\sum_{R\in \widetilde \fR(\ka,k)} \mu^{i,n}_k*\sum\Sb r(w)=r\\ R_{\ka,k}(w)=R\endSb 
b^{n,\ka}_w\Big\|_2^2 \lc 
n^{2+3(d+3)/\eta} 2^{-(k-\ka)c_0} \alpha
\sum\Sb r(w)=r\endSb 
\|b^{n,\ka}_w\|_1
\tag 5.22
$$
where $$c_0= \frac{a}{2}\min\big\{1, \frac{\tau-A}{A}\big\}
\tag 5.23$$

Given (5.22), the proof of (5.18) is a quick consequence.
First note that  $\fR(\ka, k)$ can be split into $O(10^d)$ families 
of type $\widetilde \fR(\ka,k)$. Thus 
  Minkowski's inequality and (5.22) imply that
(5.22) holds also with $\widetilde\fR(\ka, k)$ replaced 
by $\fR(\ka, k)$.
Then we obtain
from (5.20) and the modified  (5.22) that the left  side of  
(5.18) is controlled by

$$
\align &\sum\Sb n,l\endSb\sum_{k\in (I^{n}_{l})^*} \sum\Sb \ka\in (I^n_l)^*:\\\ka<
k-C_0\log n\endSb
\sum_{i\in \cJ_1^n}
\sum\Sb r\in I^{n,\ka}_l\endSb
n^{6(1+\frac{d-1}{\eta})}
 2^{-(k-\ka)c_0} \alpha
\sum\Sb r(w)=r\endSb 
\|b^{n,\ka}_w\|_1
\\
&\lc\sum\Sb n,l,\ka\endSb 
n^{6+\frac{9(d-1)}{\eta}} 
\sum\Sb k\ge \ka+C_0\log n\endSb
 2^{-(k-\ka)c_0} \alpha
\sum\Sb r\in I^{n}_l\\r<\ka\endSb\sum\Sb r(w)=r\endSb 
\|b^{n,\ka}_w\|_1.
\endalign
$$
Now we sum the geometric series
$$\sum\Sb k\ge \ka+C_0\log n\endSb
 2^{-(k-\ka)c_0} \lc n^{-c_0 C_0}$$ and using (5.23) and our choice of $C_0$ in (5.10)
we observe that $n^{-c_0 C_0}\le n^{-50(1+d/\eta)}$; this yields
that the left side of (5.18) is controlled by
$$
\alpha\sum\Sb n,l,\ka\endSb 
\sum\Sb r\in I^{n,\ka}_l\endSb\sum\Sb r(w)=r\endSb 
\|b^{n,\ka}_w\|_1 \lc \alpha\|f\|_1.$$
Thus the proof will be finished when inequality (5.22) is verified.

\subheading{Proof of (5.22)}

We   split  for fixed $n,l$, $k,\ka\in (I^n_l)^*$, $i\in \cJ^n_1$ and
 $r\in I^{n,\ka}_l$,
$$
\Big\|\sum_{R\in \widetilde \fR(\ka,k)} \mu^{i,n}_k*\sum\Sb r(w)=r\\ 
R(w)=R\endSb 
b^{n,\ka}_w\Big\|_2^2 = I+II
$$
where
$$
\align
I&= 
\sum_{R\in \widetilde \fR(\ka,k)} \int
\widetilde{ \mu^{i,n}_k}*\mu^{i,n}_k*
\sum\Sb r(w)=r\\R(w)=R\endSb 
b^{n,\ka}_w(x)   \sum\Sb r(w')=r\\R(w')=R\endSb 
\overline{b^{n,\ka}_{w'}(x)}\, dx \tag 5.24
\\
 II&= 
\sum\Sb R,R'\in \widetilde \fR(\ka,k)\\R\neq R'\endSb \int
\widetilde{ \mu^{i,n}_k}*\mu^{i,n}_k*
\sum\Sb r(w)=r\\R(w)=R\endSb 
b^{n,\ka}_w(x)   
\sum\Sb r(w')=r\\R(w')=R'\endSb 
\overline{b^{n,\ka}_{w'}(x)} dx.
\tag 5.25
\endalign
$$

We shall first estimate $II$. Fix $w$, $w'$ occuring in the expression 
(5.25).
Then using the cancellation of the $b^{n,\ka}_w$ we get
$$
\align
&\big|
\widetilde{\mu^{i,n}_k}*\mu^{i,n}_k
*b^{n,\ka}_w(x)   \big|
\\&=\Big|\int 2^{-k\tau}\Big[
\widetilde{ \mu^{i,n}_0}*\mu^{i,n}_0(\delta_{-k}(x-y))-
\sum_{j=0}^{d-1} \frac{1}{j!}
\inn{\delta_{-k}(y-x_w)}{\nabla}^j
\widetilde{\mu^{i,n}_0}*\mu^{i,n}_0 (\delta_{-k}(x-x_w))\Big]
 b^{n,\ka}_w(y)\,dy
\\&=\Big|\int_0^1 \frac{(1-s)^{d-1}}{(d-1)!} \int 2^{-k\tau}
\inn{\delta_{-k}(y-x_w)}{\nabla}^d
\widetilde{\mu^{i,n}_0}*\mu^{i,n}_0 (\delta_{-k}(x-x_w+sx_w-sy))\Big]
 b^{n,\ka}_w(y)\,dy ds
\Big|
\\&\lc 
n^{\frac 3\eta( 2d-(d-3))}
	\int_0^1
\int_w\frac{2^{-k\tau}|\delta_{-k}(y-x_w)|^d}
{|\delta_{-k}(x-x_w+sx_w-sy)|^{d+1}} 
|b^{n,\ka}_w(y)| dy ds \tag 5.26
\endalign
$$
by Lemma 2.5 applied to the measure $\mu^{i,n}_0$, with $\beta= n^{-3/\eta}$.

Now if $x\in w'$, $y\in w$ with $w'\cap R'\neq\emptyset$, 
$w\cap R\neq\emptyset$, and if 
$R\neq R'$ 
 then by the separation property of the sets in $\widetilde \fR(\ka, k)$
$$|\delta_{-k}(x-x_w)|\ge 
c_1 2^{(M_{\ka, k}-k)A} 
|\delta_{-M_{\ka, k}}(x-x_w)|\ge  10 c_1 2^{(M_{\ka, k}-k)A} 
\ge 5C_1 2^{-(k-\ka)a/2} \tag 5.27
$$
while 
$$
|\delta_{-k}(y-x_w)|
+
|\delta_{-k}(x-x_{w'})| \le 2C_1 2^{-a(k-r)}
\le 2C_1 2^{-a(k-\ka)}.
$$
Thus for $x\in w'$  we may replace
$|\delta_{-k}(x-x_w+sx_w-sy)|$ in the denominator of (5.26) by 
$|\delta_{-k}(x_w'-x_w)|$.
We also take into account that $\|b^{n,\ka}_{w'}\|_1\lc\alpha |w'|$ and thus
obtain the bound

$$
II\lc n^{3(d+3)/\eta}
\sum_{R\in \widetilde \fR(\ka,k)}\sum\Sb R(w)=R\\r(w)=r\endSb\|b^{n,\ka}_{w}\|_1
\sum\Sb R'\in \widetilde \fR(\ka,k)\\ R\neq R'\endSb
\sum\Sb R(w')=R'\\r(w')=r\endSb
 \alpha|w'| \frac{2^{-k\tau} 2^{-(k-r)ad}}
{|\delta_{-k}(x_{w'}-x_w)|^{d+1}}.
\tag 5.28
$$
Now 
we calculate using (5.27)
$$\align
&\sum\Sb R'\in \widetilde \fR(\ka,k)\\ R\neq R'\endSb
\sum\Sb R(w')=R'\\r(w')=r\endSb |w'| \frac{2^{-k\tau} 2^{-(k-r)ad}}
{|\delta_{-k}(x_{w'}-x_w)|^{d+1}}
\lc 2^{-(k-r)ad}
 \sum\Sb R'\in \widetilde \fR(\ka,k)\\ R\neq R'\endSb
\sum\Sb R(w')=R'\\r(w')=r\endSb \int_{\delta_{-k}(-x_w+w')} 
|u|^{-d-1} du
\\
&\lc 2^{-(k-\ka)ad} \int_{|u|\ge 2^{-(k-\ka)a/2}}|u|^{-d-1} du
\lc 2^{-\frac {a}{2}(k-\ka)(2d-1)}.
\endalign
$$
Combining this with (5.28) yields the bound
$$II\lc n^{3\frac{d+3}{\eta}}
 2^{-\frac {a}{2}(k-\ka)(2d-1)} \sum_{R\in \widetilde \fR(\ka, k)}\sum\Sb R(w)=R\\r(w)=r\endSb\|b^{n,\ka}_w\|_1
\tag 5.29
$$
which is controlled by the right hand side of (5.22).

We now estimate the contribution $I$.
Unfortunately, in introducing the cancellation and passing from 
$f^{n,\ka}_w$ to $b^{n,\ka}_w$
we have obscured the geometrical information on the thickness of 
$f^{n,\ka}_w$.
As the cancellation is not needed anymore for $I$ we (partially)  undo
 it and
 estimate
$$
I\le I_{1}+ I_{2}
$$
where
$$\align
I_1&= 
\sum_{R\in \widetilde \fR(\ka,k)} \int
\widetilde{ \mu^{i,n}_k}*\mu^{i,n}_k
*\sum\Sb r(w)=r\\R(w)=R\endSb 
f^{n,\ka}_w(x)   \sum\Sb r(w')=r\\R(w')=R\endSb 
\overline{b^{n,\ka}_{w'}(x)}\, dx \tag 5.30
\\
I_2&= 
\sum_{R\in \widetilde \fR(\ka,k)} \int
\widetilde{ \mu^{i,n}_k}*\mu^{i,n}_k
*\sum\Sb r(w)=r\\R(w)=R\endSb 
g^{n,\ka}_w(x)   \sum\Sb r(w')=r\\R(w')=R\endSb 
\overline{b^{n,\ka}_{w'}(x)}\, dx. \tag 5.31
\endalign
$$
Since $|g^{n,\ka}_w(x)|\lc\alpha \chi_w(x)$ we get
$$
|I_2|\lc  n^{-3\frac{d-3}{\eta}}
\alpha
\sum_{R\in \widetilde \fR(\ka,k)}  \sum\Sb r(w)=r\\R(w)=R\endSb
  \sum\Sb r(w')=r\\R(w')=R\endSb 
\int_{w'} 
\int_w\frac{2^{-k\tau}}{|\delta_{-k}(x-y)|}  dy
|b^{n,\ka}_{w'}(x)|\, dx \tag 5.32
$$
and
$$
\align
&\sum\Sb r(w)=r\\R(w)=R\endSb 
\int_w\frac{2^{-k\tau}}{|\delta_{-k}(x-y)|}  dy
\lc  2^{-(k-M(\ka,k))(\tau-A)}
\int_R\frac{2^{-M(\ka,k)\tau}}{|\delta_{-M(\ka,k)}(x-y)|}  dy
\\
&\lc  2^{-(k-\ka)(\tau-A)\frac{a}{2A}} \int_{|u|\lc 1} |u|^{-1} du
\lc  2^{-(k-\ka)(\tau-A)\frac{a}{2A}}.
\endalign
$$
Thus
$$\align
|I_2|&\lc n^{-3\frac{d-3}{\eta}} \alpha 2^{-(k-\ka)(\tau-A)\frac{a}{2A}}
\sum_{R\in \widetilde \fR(\ka,k)}  
  \sum\Sb r(w')=r\\R(w')=R\endSb \|b^{n,\ka}_{w'}\|_1.
\tag 5.33
\endalign
$$
Finally for the main term $I_1$ we use 
Lemma 2.6, then (2.28) and then  part (iii) of Proposition 4.1 to bound 
$$\align
&\Big| 
\widetilde{ \mu^{i,n}_k}*\mu^{i,n}_k
*\sum\Sb r(w)=r\\R(w)=R\endSb 
f^{n,\ka}_w(x)\Big|
\\
&=\Big|\int 
\widetilde{ \mu^{i,n}_0}*\mu^{i,n}_0(\delta_{-k} x-y)
\sum\Sb r(w)=r\\R(w)=R\endSb 
f^{n,\ka}_w(\delta_k y) dy\Big|
\\
&\lc n^{1-3(d-3)/\eta}
\Theta_n\big[
\sum\Sb r(w)=r\\R(w)=R\endSb 
f^{n,\ka}_w(\delta_k \cdot)\big]
\\
&\lc n^{1-3(d-3)/\eta} 2^{-(k-M(\ka,k))(\tau-A)}\Theta_n\big[
\sum\Sb r(w)=r\\R(w)=R\endSb 
f^{n,\ka}_w(\delta_{M(\ka, k)} \cdot)\big]
\\
&\lc
n^{2-3(d-3)/\eta} 2^{-(k-M(\ka,k))(\tau-A)} \alpha.
\endalign
$$
Since
$k-M(\ka,k)\ge(k-\ka)a/2A$ we obtain
$$
|I_1|\lc \alpha
 2^{-(k-\ka)(\tau-A)\frac{a}{2A}}
n^{2-3\frac{d-3}{\eta}}\sum_{R\in \widetilde \fR(\ka, k)}
\sum\Sb r(w')=r\\R(w')=R\endSb 
\|b^{n,\ka}_{w'}\|_1.\tag 5.34
$$
(5.33/34) and (5.29) certainly imply (5.22). This concludes the proof of Theorem 1.1.\qed

\remark{\bf Remark}The above argument also applies to maximal functions associated 
to certain surfaces with low codimension, for example if  we assume 
that for every normal vector the Gaussian 
curvature is bounded away from zero. In this case we have to work with the notions $\Lambda_{n, \beta}$, $\Theta_{n,\beta}$ in 
the remark following the proof of Proposition 2.1; here $\beta$ 
is the codimension. 
The condition  about nonvanishing  Gaussian curvature is never satisfied for 
manifolds with high codimension such as curves in three or more 
dimensions.
In those cases it is presently open whether the
weak type $L\log L$ inequality of Corollary 3.1 above can be improved.
\endremark

\head{\bf 6. Estimates for the singular integral operators}\endhead

The proof of the  weak type $L\log\log  L$  estimate 
for the singular Radon 
transforms relies to a large extent on the same arguments as 
for the maximal operator. We shall  just indicate the 
necessary modifications.

We need to prove inequality (3.21). The definition of the
 exceptional set $\cO$ and estimate (5.5) 
remains the same. Thus we are left to show
 (again with $\Phi(s)=s\log\log(e^2+s)$)
$$\meas\Big(\Big\{x:\Big|
\sum\Sb n,l\endSb \sum_{k\in (I^n_l)^*}  \mu_k^n*
\sum\Sb \ka\in (I^n_l)^*\\ \ka\le k\endSb\sum_{r(w)\in I^{n,\ka}_l}
f^{n,\ka}_w
\Big|
>\frac45 \alpha \Big\}\Big) \lc \int\Phi\big(\frac{|f(x)|}{\alpha}\big) dx.
\tag 6.1
$$

Now, as in \S5, we wish to decompose the measure 
into a part 
with curvature and a part with flatness (with the splitting depending on $n$). Some care is needed now since we need  to preserve the
cancellation of the measure when acting 
on the  $\alpha$-bounded contributions.
Before doing this decomposition we shall reverse the order of 
the steps (5.9), (5.11) and first 
get  an  analogue of (5.11) for the functions
$\mu_k^n$. Indeed since $\|\mu^n_k\|_1=O(1)$ the  argument for (5.11) 
yields
$$\meas\Big(\Big\{x:\Big|
\sum\Sb n,l,\ka\endSb\sum\Sb k\in (I^n_l)^*\\\ka\le k\le \ka+C_0\log n
\endSb
\mu^n_k*\sum\Sb r(w)\in I^{n,\ka}_l\endSb f^{n,\ka}_w\Big|>
\frac{\alpha}{10}\Big\}\Big)
\lc 
\int\Phi\big(\frac{|f(x)|}{\alpha}\big) dx
\tag  6.2
$$
and therefore we have to bound
$$\meas\Big(\Big\{x:\Big|
\sum\Sb n,l\endSb
\sum\Sb \ka\in (I^n_l)^*\endSb
\sum\Sb k\in (I^n_l)^*\\ k\ge \ka+C_0\log n
\endSb
\mu^n_k*\sum\Sb r(w)\in I^{n,\ka}_l\endSb 
f^{n,\ka}_w
\Big|>\frac 7{10}\alpha\Big\}\Big).
\tag  6.3
$$

As before we split 
$f^{n,\ka}_w=g^{n,\ka}_w+b^{n,\ka}_w$
and we first show that
$$\meas\Big(\Big\{x:\Big|
\sum\Sb n,l\endSb
\sum\Sb \ka\in (I^n_l)^*\endSb
\sum\Sb k\in (I^n_l)^*\\ k\ge \ka+C_0\log n
\endSb
\mu^n_k*\sum\Sb r(w)\in I^{n,\ka}_l\endSb 
g^{n,\ka}_w
\Big|>\frac{\alpha}{10}\Big\}\Big) \lc
\frac{\|f\|_1}{\alpha}.
\tag  6.4
$$

We use the nonisotropic version of  an inequality  in \cite{6, p. 548} for the maximal version of the singular integral, namely we have
$$
\Big\|\sup_{K_1,K_2} \Big|\sum_{k=K_1}^{K_2}\widetilde 
 \mu_k*u\Big|\Big\|_2\lc
\|u\|_2. \tag 6.5
$$
Here
$\widetilde  \mu_k$ is the reflection of 
$ \mu_k$.
 Indeed 
for (6.5)  one just  needs
$|\widehat \mu(\xi)|\le \min\{|\xi|, |\xi|^{-\gamma}\}$
for some $\gamma>0$ ({\it cf.} (3.15)).
In order to use (6.5) we have to split 
$\mu_k^n = \mu_k-(\mu_k-\mu_k^n)$.

From (6.5) and (5.17) we get
$$\align
&\Big\|\sum\Sb n,l\endSb
\sum\Sb \ka\in (I^n_l)^*\endSb
\sum\Sb k\in (I^n_l)^*\\ k\ge \ka+C_0\log n
\endSb
\mu_k*\sum\Sb r(w)\in I^{n,\ka}_l\endSb 
g^{n,\ka}_w\Big\|_2^2
\\&=\sup\limits_{\|u\|_2\le 1}\Big|
\int 
\sum\Sb n,l\endSb
\sum\Sb \ka\in (I^n_l)^*\endSb 
\sum\Sb r(w)\in I^{n,\ka}_l\endSb 
g^{n,\ka}_w(y)
\sum\Sb k\in (I^n_l)^*\\ k\ge \ka+C_0\log n
\endSb 
\widetilde \mu_k*u(y) dy\Big|^2
\\&\le
\sup\limits_{\|u\|_2\le 1}\Big\|
\sup_{K_1,K_2}\Big|\sum_{k=K_1}^{K_2}
\widetilde \mu_k*u\Big|\Big\|_2^2
\Big\|
\sum\Sb n,l\endSb
\sum\Sb \ka\in (I^n_l)^*\endSb 
\Big|\sum\Sb r(w)\in I^{n,\ka}_l\endSb 
g^{n,\ka}_w\Big|\Big\|_2^2 
\\
&\lc\Big\|
\sum\Sb n,l\endSb
\sum\Sb \ka\in (I^n_l)^*\endSb 
\sum\Sb r(w)\in I^{n,\ka}_l\endSb 
|g^{n,\ka}_w|\Big\|_2^2.
\tag 6.6
\endalign
$$
For each $w$ and $x\in w$ we have 
$$
\sum\Sb n,l\endSb
\sum\Sb \ka\in (I^n_l)^*\\ r(w)\in I^{n,\ka}_l\endSb
|g^{n,\ka}_w(x)|\lc
\sum\Sb n,l\endSb
\sum\Sb \ka\in (I^n_l)^*\\ r(w)\in I^{n,\ka}_l\endSb \frac{1}{|w|}\int_w
|f^{n,\ka}_w(x)|dx \lc \frac{1}{|w|}\int_w |f(x)| dx\lc \alpha
$$
and in view of the disjointness of the sets $w$ the expression (6.6) is 
controlled by
$$\sum_w\Big\|
\sum\Sb n,l\endSb
\sum\Sb \ka\in (I^n_l)^*\endSb 
\sum\Sb r(w)\in I^{n,\ka}_l\endSb 
|g^{n,\ka}_w|\Big\|_2^2\lc \alpha \|f\|_1.
\tag 6.7
$$
Moreover for fixed $n$, and $m\ge n$ we get using (3.15)
$$
\align
&\Big\|
\sum\Sb l\endSb
\sum\Sb k\in (I^n_l)^*\endSb
(\mu_k^m-\mu_k^{m+1})*
\sum\Sb \ka\in (I^n_l)^*\endSb 
\sum\Sb r(w)\in I^{n,\ka}_l\endSb 
g^{n,\ka}_w\Big\|_2 
\\
&\lc 
2^{-m\ga}
\sum\Sb l\endSb
\sum\Sb k\in (I^n_l)^*\endSb
\Big\|
\sum\Sb \ka\in (I^n_l)^*\endSb 
\sum\Sb r(w)\in I^{n,\ka}_l\endSb 
g^{n,\ka}_w\Big\|_2^2 
\endalign
$$
and thus using the telescoping sum 
$\mu_k^n-\mu_k=\sum_{m=n}^\infty
(\mu_k^m-\mu_k^{m+1})$ we obtain
$$
\align
&\Big\|
\sum\Sb n, l\endSb
\sum\Sb k\in (I^n_l)^*\endSb
(\mu_k^n-\mu_k)*
\sum\Sb \ka\in (I^n_l)^*\endSb 
\sum\Sb r(w)\in I^{n,\ka}_l\endSb 
g^{n,\ka}_w\Big\|_2 
\\
&\lc
\Big[\sum_n\sum_{m=n}^{\infty} 2^{-n\gamma}
\Big(
\sum\Sb l\endSb
\sum\Sb k\in (I^n_l)^*\endSb
\Big\|
\sum\Sb \ka\in (I^n_l)^*\endSb 
\sum\Sb r(w)\in I^{n,\ka}_l\endSb 
g^{n,\ka}_w\Big\|_2^2 \Big)^{1/2}\Big]
\\
&\lc\Big[\sum_n 2^{-n\gamma}n^2
\Big(
\sum\Sb l\endSb
\sum\Sb w: r(w)\\ \in I^n_l\endSb
\Big\|
\sum\Sb \ka\in (I^n_l)^*\endSb 
g^{n,\ka}_w\Big\|_2^2 \Big)^{1/2}\Big]
\endalign
$$
which by the argument above is dominated by a constant times 
$\alpha\|f\|_1$.
Combining these  estimates with Chebyshev's inequality 
 we see that (6.4) holds.

We are left to prove 
$$\meas\Big(\Big\{x:\Big|
\sum\Sb n,l\endSb
\sum\Sb \ka\in (I^n_l)^*\endSb
\sum\Sb k\in (I^n_l)^*\\ k\ge \ka+C_0\log n
\endSb
\mu^n_k*\sum\Sb r(w)\in I^{n,\ka}_l\endSb 
b^{n,\ka}_w
\Big|>\frac{\alpha}{10}\Big\}\Big) \lc
\frac{\|f\|_1}{\alpha}.
\tag  6.8
$$

We now let
$\mu^{i,n,m}=\nu^{i,n}*\phi_m$ 
(which was previously considered only for the case $m=n$)
and define the $L^1$ dilate
$\mu^{i,n,m}_k=2^{-k\tau}\mu^{i,n,m}(\delta_{-k}\cdot)$.
Split (with $\cJ^n_1$ and $\cJ^n_2$ as in  \S5)
$$\align\mu^n_k&=\mu^0_k+(\mu^n_k-\mu^{0}_k)
\\
&= \mu^0_k+\sum_{i\in \cJ^n_1}\sum_{m=1}^n(\mu^{i,n,m}_k-\mu^{i,n,m-1}_k)
+\sum_{i\in \cJ^n_2}(\mu^{i,n,n}_k-\mu^{i,n,0}_k).
\tag 6.9
\endalign
$$ 
Let
$h^{i,n,m}_k=\mu^{i,n,m}_k-\mu^{i,n,m-1}_k$ and
$h^{n,m}_k=\sum_{i\in \cJ_1^n} h^{i,n,m}_k$. 
Using (6.9) we  split 
$$
\sum\Sb n,l\endSb
\sum\Sb \ka\in (I^n_l)^*\endSb
\sum\Sb k\in (I^n_l)^*\\ k\ge \ka+C_0\log n
\endSb
\mu^n_k*\sum\Sb r(w)\in I^{n,\ka}_l\endSb b^{n,\ka}_w= I+II+\sum_{m=1}^\infty III_m
$$
where
$$\align
I&=\sum\Sb n,l\endSb
\sum\Sb \ka\in (I^n_l)^*\endSb
\sum\Sb k\in (I^n_l)^*\\ k\ge \ka+C_0\log n\endSb
\mu^0_k*\sum\Sb r(w)\in I^{n,\ka}_l\endSb b^{n,\ka}_w
\\
II&=
\sum\Sb n,l\endSb
\sum\Sb \ka\in (I^n_l)^*\endSb
\sum\Sb k\in (I^n_l)^*\\ k\ge \ka+C_0\log n\endSb
\sum_{i\in \cJ^n_2}
(\mu^{i,n,n}_k-\mu^{i,n,0}_k)
*\sum\Sb r(w)\in I^{n,\ka}_l\endSb b^{n,\ka}_w
\\
III_m
&=
\sum\Sb n\ge m\endSb\sum_l
\sum\Sb \ka\in (I^n_l)^*\endSb
\sum\Sb k\in (I^n_l)^*\\ k\ge \ka+C_0\log n\endSb
\sum_{i\in \cJ^n_1}
(\mu^{i,n,m}_k-\mu^{i,n,m-1}_k)
*\sum\Sb r(w)\in I^{n,\ka}_l\endSb b^{n,\ka}_w.
\endalign
$$

We show that 
$$
 \|I\|_{L^1(\bbR^d\setminus \Omega^*)}
+\|II\|_1
\lc 
\|f\|_1
\tag 6.10
$$
and
$$\align
&\|III_m\|_2^2
\lc (1+m)^{-2} \alpha\|f\|_1.
\tag 6.11
\endalign
$$
(6.10/11) imply that the sets where
$|I|>\alpha/10$,
$|II|>\alpha/10$,
and $\sum_{m=1}^\infty| III_m|>\alpha/10$ all have measure 
$\lc \alpha^{-1}\|f\|_1$. Combining this with the estimate (3.2)
 for the measure of $\Omega^*$  yields (6.8).

The inequality  $$ \|I\|_{L^1(\bbR^d\setminus \Omega^*)}
\lc  \|f\|_1$$
follows from the standard estimates for singular integrals
 (in view of the 
regularity of $\mu*\phi_0$). The bound for $\|II\|_1$ is 
proved exactly as in estimate (5.9).
Thus we are left to check (6.11).

Concerning the terms $III_m$ we apply Cauchy-Schwarz' inequality and 
estimate
$$
\|III_m\|_2^2\lc 
\sum\Sb n\ge m\endSb n^2\Big\| \sum_l
\sum\Sb \ka\in (I^n_l)^*\endSb
\sum\Sb k\in (I^n_l)^*\\ k\ge \ka+C_0\log n\endSb
\sum_{i\in \cJ^n_1}
h^{i,n,m}_k
*\sum\Sb r(w)\in I^{n,\ka}_l\endSb b^{n,\ka}_w\Big\| \lc C_2IV_m+V_m
$$
where $C_2>10/a$,
$$IV_m=
\sum\Sb n\ge m\endSb  n^2\sum_l\Big\| \sum\Sb \ka\in (I^n_l)^*\endSb
\sum\Sb k\in (I^n_l)^*\\ k\ge \ka+C_0\log n\endSb
\sum_{i\in \cJ^n_1}
h^{i,n,m}_k
*\sum\Sb r(w)\in I^{n,\ka}_l\endSb b^{n,\ka}_w\Big\|_2^2 
\tag 6.12
$$
and $$
V_m=
\sum\Sb n\ge m\endSb  n^2\sum\Sb l,l'\\|l-l'|\ge  C_2\endSb
\sum\Sb \ka,\ka'\in (I^n_l)^*\endSb
\sum\Sb k,k'\in (I^n_l)^*\\ k\ge \ka+C_0\log n  
\\k'\ge \ka'+C_0\log n\endSb
\sum\Sb i, i'\in \cJ^n_1\endSb 
\biginn{
\widetilde h^{i',n,m'}_{k'} *h^{i,n,m}_k
*\sum\Sb r(w)\in I^{n,\ka}_l\endSb b^{n,\ka}_w}
{\sum\Sb r(w')\in I^{n,\ka'}_{l'}\endSb  \overline{ b^{n,\ka'}_{w'}}}.
\tag 6.13
$$

The inner product in the
second term   is estimated by  Plancherel's theorem.
By van der Corput's Lemma and cancellation
there is the  Fourier transform estimate
$$\big|\Cal F[h_0^{i,n,m}](\xi)\big|\lc\min\{|\xi|,|\xi|^{-\gamma}\}$$
and thus
$$\big|\Cal F\big[\widetilde h^{n,m, i'}_{k'}*h^{n,m, i}_k]\big|\lc 
2^{-|k-k'|a\gamma}$$
which is   $O(2^{-n|l-l'|\gamma a/2})$ if 
$k\in (I^n_l)^*$, $k'\in (I^n_{l'})^*$ , $|l-l'|\ge C_2>10/a$. 
Set
$$\cE_{l,n}(x)=
\sum\Sb \ka\in (I^n_l)^*\endSb\sum\Sb r(w)\in I^{n,\ka}_l\endSb b^{n,\ka}_w (x).
$$
We may apply Cauchy-Schwarz and Parseval's theorem to
bound
$$
\align
V_m&\lc
\sum\Sb n\ge m\endSb n^{4+6(d-1)/\eta}  \sum\Sb l,l'\\|l-l'|\ge C_2\endSb 
2^{-n|l-l'|\gamma}\int|
\widehat{\cE_{l,n}}(\xi)||
\widehat{\cE_{l',n}}(\xi)|d\xi
\\
&\lc
\sum\Sb n\ge m\endSb n^{4+6(d-1)/\eta}  \sum\Sb l,l'\\|l-l'|\ge C_2\endSb 
2^{-n|l-l'|\gamma}\|\cE_{l,n}\|_2\|\cE_{l',n}\|_2
\\&\lc
\sum\Sb n\ge m\endSb n^{4+6(d-1)/\eta}  2^{-C_2 a\gamma n}\sum_l
\|\cE_{l,n}\|_2^2
\endalign
$$
Now 
$$
\|\cE_{l,n}\|_2^2\lc 2^{c(n+1)}\alpha 
\Big\|
\sum\Sb r(w)\in I^{n,\ka}_l\endSb |b^{n,\ka}_w|\Big\|_1
$$
where $c$ is as in (3.3) and hence we obtain
$$
V_m\lc
\alpha\sum\Sb n\ge m\endSb 2^{-5 a \gamma n}
\Big\|\sum_l\sum\Sb \ka\in (I^n_l)^*\endSb
\sum\Sb r(w)\in I^{n,\ka}_l\endSb |b^{n,\ka}_w|\Big \|_1\lc 
2^{-5a\gamma m} \alpha\|f\|_1.
\tag 6.14
$$

For the  term $IV_m$  
we have by Cauchy-Schwarz for the $k$ summation and other applications 
of Cauchy-Schwarz leading to (5.20)
$$
\align
IV_m&\lc
\sum\Sb n\ge m\endSb  n^3\sum_l\sum\Sb k\in (I^n_l)^*\endSb
\Big\| \sum\Sb \ka\in (I^n_l)^*\\ k\ge \ka+C_0\log n\endSb
\sum_{i\in \cJ^n_1}
h^{i,n,m}_k
*\sum\Sb r(w)\in I^{n,\ka}_l\endSb b^{n,\ka}_w\Big\|_2^2 
\\
&\lc
\sum\Sb n\ge m\endSb  n^{5+\frac{3(d-1)}{\eta}}\sum_l
\sum\Sb k\in (I^n_l)^*\endSb
 \sum\Sb \ka\in (I^n_l)^*\\ k\ge \ka+C_0\log n\endSb
\sum_{i\in \cJ^n_1}\sum_{r\in I^{n,l}_\ka}
\Big\|
\mu^{i,n,m}_k
*\sum\Sb r(w)=r\endSb b^{n,\ka}_w\Big\|_2^2 .
\tag 6.15
\endalign
$$

Now $\mu^{i,n,m}_k$ satisfies similar quantitative properties
as $\mu^{i,n,n}_k\equiv\mu^{i,n}_k$ considered in \S5; in particular we have
$|\partial^\alpha (\mu^{i,n,n}_0 * \mu^{i,n,m}_0)(x)|\lc n^{-3(d-3-2|\alpha|)/\eta} 
(2^{-m}+|x|)^{-1-|\alpha|}$. Thus the estimates 
 for expression  (5.20), are applicable and we obtain the bound
$$\align
IV_m&\lc \alpha
\sum\Sb n\ge m\endSb  n^{8+9(d+1)/\eta} 
\sum_{\ka} \sum\Sb k\ge \ka+C_0\log n\endSb 2^{-(k-\ka)c_0}
\sum_r \sum\Sb r(w)=r\endSb \|b^{n,\ka}_w\|_1
\\&\lc \sum_{n\ge m} n^{8+9(d+1)/\eta -C_0c_0}  \alpha\|f\|_1 \lc (1+m)^{-2}\alpha\|f\|_1.
\tag 6.16
\endalign
$$
This shows (6.11) and thus (6.8) and the proof of   Theorem 1.2. is complete.\qed

\enddocument